\documentclass[twocolumn]{IEEEtran}
\usepackage{mathtools,amsmath,amssymb,amsthm,epsfig,color,empheq,graphicx}
\usepackage{enumerate,url,algpseudocode,algorithm,enumitem,gensymb}
\usepackage{wasysym,epstopdf}
\usepackage[caption=false,font=footnotesize]{subfig}
\usepackage[table]{xcolor}
\usepackage{arydshln,tikz,balance}
\usetikzlibrary{matrix,decorations.pathreplacing}
\usepackage{dashbox, hyperref}

\DeclareMathOperator{\rank}{rank}

\DeclareMathOperator{\trace}{Tr}

\newtheorem{lemma}{Lemma}

\newtheorem{remark}{Remark}

\newcommand \bzero{\mathbf{0}}
\newcommand \bone{\mathbf{1}}
\newcommand \ba{\mathbf{a}}
\newcommand \bb{\mathbf{b}}
\newcommand \bc{\mathbf{c}}
\newcommand \bd{\mathbf{d}}
\newcommand \be{\mathbf{e}}

\newcommand \bp{\mathbf{p}}
\newcommand \bq{\mathbf{q}}

\newcommand \bs{\mathbf{s}}

\newcommand \bu{\mathbf{u}}
\newcommand \bv{\mathbf{v}}

\newcommand \bx{\mathbf{x}}
\newcommand \by{\mathbf{y}}

\newcommand \bA{\mathbf{A}}

\newcommand \bC{\mathbf{C}}
\newcommand \bD{\mathbf{D}}
\newcommand \bE{\mathbf{E}}

\newcommand \bG{\mathbf{G}}

\newcommand \bJ{\mathbf{J}}
\newcommand \bK{\mathbf{K}}
\newcommand \bL{\mathbf{L}}
\newcommand \bM{\mathbf{M}}
\newcommand \bN{\mathbf{N}}

\newcommand \bV{\mathbf{V}}

\newcommand \bY{\mathbf{Y}}

\newcommand \btheta{\boldsymbol{\theta}}

\newcommand \mcA{\mathcal{A}}

\newcommand \mcG{\mathcal{G}}

\newcommand \mcM{\mathcal{M}}
\newcommand \mcN{\mathcal{N}}
\newcommand \mcO{\mathcal{O}}
\newcommand \mcP{\mathcal{P}}

\newcommand \mcS{\mathcal{S}}
\newcommand \mcT{\mathcal{T}}

\newcommand \tbv{\tilde{\mathbf{v}}}

\newcommand \tbx{\tilde{\mathbf{x}}}

\newcommand \tbV{\tilde{\mathbf{V}}}

\newcommand \tbY{\tilde{\mathbf{Y}}}

\newcommand \hbv{\hat{\mathbf{v}}}

\newcommand \hbx{\hat{\mathbf{x}}}

\usepackage{xr-hyper}
\allowdisplaybreaks

\begin{document}
\title{Smart Inverter Grid Probing for Learning Loads:\\
Part II -- Probing Injection Design}

\author{
	Siddharth Bhela,~\IEEEmembership{Student Member,~IEEE,}
	Vassilis Kekatos,~\IEEEmembership{Senior Member,~IEEE,} and
	Sriharsha Veeramachaneni
	
%
}


\maketitle
\begin{abstract}
This two-part work puts forth the idea of engaging power electronics to probe an electric grid to infer non-metered loads. Probing can be accomplished by commanding inverters to perturb their power injections and record the induced voltage response. Once a probing setup is deemed topologically observable by the tests of Part I, Part II provides a methodology for designing probing injections abiding by inverter and network constraints to improve load estimates. The task is challenging since system estimates depend on both probing injections and unknown loads in an implicit nonlinear fashion. The methodology first constructs a library of candidate probing vectors by sampling over the feasible set of inverter injections. Leveraging a linearized grid model and a robust approach, the candidate probing vectors violating voltage constraints for any anticipated load value are subsequently rejected. Among the qualified candidates, the design finally identifies the probing vectors yielding the most diverse system states. The probing task under noisy phasor and non-phasor data is tackled using a semidefinite-program (SDP) relaxation. Numerical tests using synthetic and real-world data on a benchmark feeder validate the conditions of Part I; the SDP-based solver; the importance of probing design; and the effects of probing duration and noise.
\end{abstract}



\section{Introduction}\label{sec:introd}
Part I of this work put forth the novel data acquisition scheme of probing-to-learn (P2L). The P2L scheme leverages smart inverters to probe an electric grid with the purpose of finding the values of non-metered loads. It also provided conditions under which a particular probing setup is successful. In particular, it was shown that given the feeder graph $\mcG$, the locations of non-metered buses $\mcO$ and probing buses $\mcM$, and the number of probing actions $T$, a simple linear program could tell whether non-metered loads could be recovered or not. Assuming noiseless data, this test relied on the generic rank of the Jacobian matrix $\bJ\left(\{\bv_t\}\right)$ related to the P2L equations. It is thus a \emph{topological} rather than a \emph{numerical observability} guarantee~\cite[Ch.~4.6]{ExpConCanBook}.

Even for the standard power flow (PF) and power system state estimation (PSSE) setups, topological observability relates to the sparsity structure of the associated Jacobian matrix. This structure alone however cannot adequately capture the numerical column rank of the Jacobian matrix. There exist specification or measurement sets whose Jacobian is full column-rank in general, but becomes ill-conditioned or even singular under specific state values (including the boundaries for voltage collapse); see e.g.,~\cite[Ch.~10]{ExpConCanBook}, \cite{Yang17}. In addition, once a probing setup is deemed topologically observable, the power injections of probing inverters could be judiciously selected to improve load or state estimates. This is challenging since P2L is an implicit nonlinear identification task, and probing injections should be comply to network constraints without knowing the non-metered loads. 

The contribution of Part II is on two practical aspects of grid probing. First, a systematic approach to design probing setpoints that conform to grid safety and improve numerical accuracy is developed in Section~\ref{sec:probing}. Second, the proposed P2L task is tackled through semidefinite program (SDP)-based solvers presented in Section~\ref{sec:solvers}. The conditions of Part I along with the probing setpoint design and the solver of Part II, are numerically validated using actual residential load data from the Pecan Street project on the IEEE 34-bus benchmark feeder in Section~\ref{sec:sim}. Conclusions and current research efforts are outlined in Section~\ref{sec:conclusion}.

Adding to the notational conventions of Part I, here the symbol $\bone$ denotes the all-one vector and $\be_k$ is the $k$-th canonical vector; their dimensions would be clear from the context. The notation $\bV \succeq \bzero$ means that $\bV$ is a Hermitian (complex and conjugate symmetric) positive semidefinite matrix; the matrix trace is denoted by $\trace(\cdot)$; and $\|\ba\|_2$ is the $\ell_2$-norm of vector $\ba$. The notation $k=1:K$ is a shorthand to $k=1,\dots, K$.

\section{Designing Probing Injections}\label{sec:probing}
Suppose probing has been deemed successful for a particular $(\mcM,\mcO)$ placement of probed and non-metered buses, i.e., the setup $(\mcM,\mcO)$ has passed the test of Algorithm~1 or 2 of Part I. The next question is how to select probing setpoints that are implementable by inverters; compliant to feeder constraints; and at the same time, improve estimation accuracy. This section deals with the design of inverter setpoints during probing interval $\mcT$ with slots $t=1,\ldots,T$, for a given $T$.

To facilitate the exposition, let us stack the power injections at all probing buses $\{(p_{n},q_{n})\}_{n\in\mcM}$ in vectors $\bp_\mcM$, $\bq_\mcM$, and $\bs_\mcM:=[\bp_\mcM^\top~\bq_\mcM^\top]^\top$. Likewise, the injections at all non-metered buses $\{(p_{n},q_{n})\}_{n\in\mcO}$ are collected in $\bp_\mcO$, $\bq_\mcO$, and $\bs_\mcO:=[\bp_\mcO^\top~\bq_\mcO^\top]^\top$. The injections at slot $t$ will be denoted by a superscript $t$.

In search of a meaningful metric to design the probing injections $\{\bs_\mcM^t\}_{t=1}^T$, one could consider the minimum mean square estimation error for non-metered loads $\bs_\mcO$ or states $\{\bv_t\}_{t=1}^T$. The former is hard to derive given the implicit estimation task involved. The latter exhibits the Cramer-Rao lower bound (CRLB) of $[\bJ^\top\left(\{\bv_t\}\right)\bJ\left(\{\bv_t\}\right)]^{-1}$; a proof for this CRLB can be obtained by adopting the result in~\cite{Gang18}. Since $\bJ\left(\{\bv_t\}\right)$ depends linearly on $\{\bv_t\}_{t=1}^T$, the CRLB depends inverse quadratically on the unknown states. 

To arrive at a practical solution, we resort to selecting probing setpoints so that the electric grid is driven to the most diverse states $\{\bv_t\}_{t=1}^T$ while abiding by inverter and feeder operational constraints. We conjecture that probing the grid to effect larger state variations across $\mcT$ would yield smaller condition numbers for $\bJ\left(\{\bv_t\}\right)$ and $\bJ^\top\left(\{\bv_t\}\right)\bJ\left(\{\bv_t\}\right)$. 


Hence, the goal is to design $\{\bs_\mcM^t\}_{t=1}^T$ that yield the most diverse system states $\{\bv_t\}_{t=1}^T$. Since the system states depend on both $\{\bs_\mcM^t\}_{t=1}^T$ and the unknown $\bs_\mcO$ in a non-linear fashion, our design adopts a linearized power flow model. The latter can be obtained by taking the first-order Taylor's series approximation of the PF equations with respect to nodal voltages expressed in polar coordinates~\cite{BoDo15}, \cite{Deka17}. Unless a reference system state is available, the linearization occurs at the flat voltage profile of $\tbv=u_0\bone+j\bzero$, and yields the so termed \emph{linearized distribution flow} (LDF) model~\cite{BW2}, \cite{Deka17}, which can be rearranged for our analysis as
\begin{equation}\label{eq:LDF}
\by:=\begin{bmatrix}
\bu-u_0\bone\\
\boldsymbol{\theta}
\end{bmatrix} 
= 
\begin{bmatrix}
\bK & \bL \\
\bM & \bN
\end{bmatrix}
\begin{bmatrix}
\bs_\mcM \\
\bs_\mcO
\end{bmatrix}.
\end{equation}
The vectors $\bu$ and $\btheta$ collect the voltage magnitudes and angles at all buses excluding the substation; and matrices $(\bK,\bL,\bM,\bN)$ depend on the bus admittance matrix $\bY$; see~\cite{Deka17}, \cite{GVV17}. Armed with a linear mapping between power injections and voltages, the design of setpoints $\{\bs_\mcM^t\}_{t=1}^T$ is accomplished next in three steps. 

\subsection{Build Library of Implementable Probing Setpoints}\label{subsec:S1}
The first step of the setpoint design builds a library $\mcS$ of $K\gg T$ \emph{candidate} injection vectors indexed by $k$
\begin{equation}\label{eq:S}
\mcS:=\{\bs_{\mcM}^k\}_{k=1}^K.
\end{equation}
The entries of each $\bs_{\mcM}^k$ should be \emph{implementable}, in the sense that each probing inverter should be able to inject the requested value of complex power.

To characterize the allowable range of inverter injections $(p_{n},q_{n})$ with $n\in\mcM$, two inverter classes are identified. The first class consists of inverters interfacing solar panels. When inverter $n$ interfaces a solar panel, its complex injection is limited by its apparent power capacity $\bar{s}_n$ as
\begin{equation}\label{eq:appc1}
p_{n}^2+q_{n}^2 \leq \bar{s}_n^2.
\end{equation}
Moreover, if the maximum active power that can be generated given the solar irradiance at the current probing period is $\bar{p}_{n}$, then its active power injection is limited by
\begin{equation}\label{eq:appc2}
0\leq p_{n}\leq \bar{p}_{n}.
\end{equation}

The second class consists of inverters interfacing energy storage units. The apparent power constraint of \eqref{eq:appc1} should still be enforced. If the power rate of energy storage unit $n$ is $\bar{p}_{n}$, the active injection from inverter $n$ should lie within
\begin{equation}\label{eq:appc3}
-\bar{p}_{n}\leq p_{n}\leq \bar{p}_{n}
\end{equation}
since the battery can be charged or discharged. Given the short duration of probing, limits on the state of charge have been ignored for simplicity.

Given the limitations for each inverter class, a candidate probing injection $\bs_\mcM^k\in\mcS$ can be constructed by sampling uniformly at random $p_n^k$ within \eqref{eq:appc2}--\eqref{eq:appc3} for all $n\in\mcM$. Upon fixing active injections, the reactive injections can be sampled again uniformly at random within $|q_n^k|\leq \sqrt{\bar{s}_n^2-(p_n^k)^2}$ to comply with \eqref{eq:appc1}. Scenarios where a single bus hosts multiple inverters belonging to the previous two or additional classes can be incorporated in the sampling process.

As explained in Remark~1 of Part I, a probing bus $n\in\mcM$ may be hosting controllable inverters and non-controllable assets (non-probing inverters and non-controllable loads). The process of sampling implementable injections through \eqref{eq:appc1}--\eqref{eq:appc3} can be repeated for all controllable inverters. The net injection from non-controllable assets is assumed to be metered; that is the case for the Pecan Street dataset~\cite{pecandata}. The complex powers injected into bus $n$ are summed up and used in the P2L. To keep the notation uncluttered, we will slightly abuse notation and denote this net injection at bus $n$ as $p_n+jq_n$.

The sampling process is repeated $K$ times to construct library $\mcS$. Although each candidate probing vector $\bs_\mcM^k\in\mcS$ can be implemented by inverters, the aggregate effect of probing injections may be violating feeder constraints. To handle this concern, we next reduce library $\mcS$ to only those probing injections abiding by feeder constraints.

\subsection{Maintaining only Network-Compliant Probing Setpoints}\label{subsec:S2}
Even though a probing action lasts for one second or two, the operator may still want to guarantee that it does not violate any feeder constraints. For example, voltage regulation standards dictate voltage magnitudes to remain within a pre-specified range as $\underline{u}\leq u_n\leq \overline{u}$ for all $n\in\mcN^+$. A probing injection vector $\bs_\mcM^k\in\mcS$ is deemed \emph{network-compliant} if the incurred voltage deviations are maintained within the allowable range $\underline{u}\bone\leq\bu -u_0\bone\leq \overline{u}\bone$ with the inequalities applied entry-wise. Thanks to \eqref{eq:LDF}, these voltage constraints can be expressed as linear inequality constraints on $\bs_\mcM^k$
\begin{equation}\label{eq:volt}
\underline{u}\bone\leq \bK\bs_\mcM^k +\bL \bs_\mcO \leq \overline{u}\bone.
\end{equation}

One cannot directly check whether $\bs_\mcM^k$ is network-compliant, since $\bs_\mcO$ is unknown. To bypass this complication, non-metered loads are assumed to lie within a known range
\begin{equation} \label{eq:load}
\underline{\bs}_{\mcO} \leq \bs_{{\mcO}} \leq \overline{\bs}_{\mcO}.
\end{equation}
The bounds $(\underline{\bs}_{\mcO},\overline{\bs}_{\mcO})$ can be derived from historical data, the confidence intervals of load forecasts, or the load estimates obtained during the previous probing period. 

Adopting a robust design, we would like to comply with the voltage constraints in \eqref{eq:volt} for all possible values of non-metered loads in \eqref{eq:load}. To do so, we leverage the next version of Farka's lemma on the containment of polytopes.

\begin{lemma}[\cite{Eaves1982}, \cite{Mangasarian2002}, \cite{Zhao17}]\label{le:farka}
The non-empty polytope $\mcP_1:=\{\bx:\bA\bx\leq \bb\}$ with $\bA\in\mathbb{R}^{M\times N}$ is contained within the polytope $\mcP_2:=\{\bx:\bC\bx\leq \bd\}$ with $\bC\in\mathbb{R}^{K\times N}$ if and only if there exists matrix $\bE\geq \bzero$ satisfying $\bE\bA=\bC$ and $\bE\bb\leq \bd$.
\end{lemma}

Based on Lemma~\ref{le:farka}, to ensure that the polytope over $\bs_\mcO$ defined in \eqref{eq:load} is contained within the polytope of \eqref{eq:volt}, we need to solve the feasibility problem
\begin{align}\label{eq:farkas}
\mathrm{find}~&~\bE\\
\mathrm{s.to}~&~\bE\geq \bzero\nonumber\\
~&~\bE\left[\begin{array}{c}
-\mathbf{I}_{2O}\\
\mathbf{I}_{2O}
\end{array}\right]= 
\left[\begin{array}{c}
-\mathbf{L}\\
\mathbf{L}
\end{array}\right]\nonumber\\
~&~\bE\left[\begin{array}{c}
-\underline{\bs}_\mcO\\
\overline{\bs}_\mcO
\end{array}\right]\leq 
\left[\begin{array}{c}
\bK\bs_\mcM^k+(u_0-\underline{u})\bone\\
-\bK\bs_\mcM^k-(u_0-\overline{u})\bone
\end{array}\right].\nonumber
\end{align}
Given $\bs_\mcM^k$, if the linear program in \eqref{eq:farkas} is feasible, the candidate vector $\bs_\mcM^k$ is deemed network-compliant and is copied to the \emph{reduced library} $\mcS_r$. Otherwise, the candidate vector is not copied to $\mcS_r$ since there exist load values within $[\underline{\bs}_{\mcO},\overline{\bs}_{\mcO}]$ that violate the voltage constraints in \eqref{eq:volt}. The test of \eqref{eq:farkas} is repeated for all $\bs_\mcM^k\in\mcS$ to get the reduced library $\mcS_r:=\{\bs_\mcM^\ell\}_{\ell=1}^L$ of $L$ candidate injection vectors with $L\leq K$.

\begin{figure}[t]
	\centering
	\includegraphics[scale=0.4]{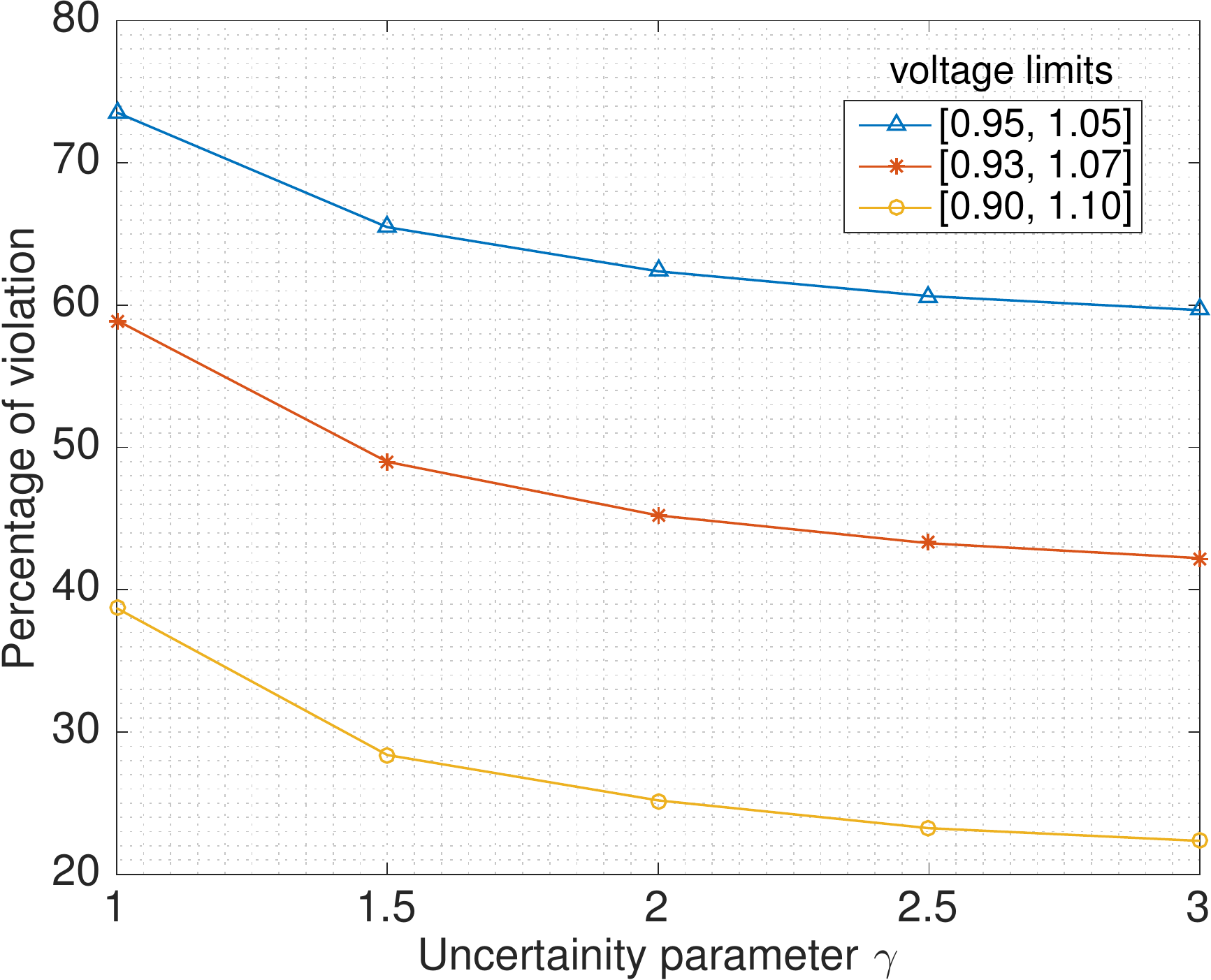}
	\caption{Percentage of candidate vectors $\bs_\mcM^k\in\mcS$ that violate \eqref{eq:volt} for varying voltage bounds $(\underline{u},\overline{u})$. The percentage of non-admissible vectors decreases with smaller load uncertainty and/or looser voltage regulation bounds.}
	\label{fig:clustT4}
\end{figure}

To demonstrate the importance of this library reduction step, we ran a numerical test on the IEEE 34-bus feeder for $T=4$ and $O=10$; see Fig.~\ref{fig:clustT4}. Load uncertainty in \eqref{eq:load} was confined within $\underline{\bs}_{\mcO}=(1-\frac{1}{\gamma})\bs_{\mcO}$ and $\overline{\bs}_{\mcO}=(1+\frac{1}{\gamma})\bs_{\mcO}$ for $\gamma>0$. The candidate inverter injections in $\mcS$ were randomly drawn from $\pm 0.2$ pu and tested against \eqref{eq:farkas}. For increasing $\gamma$, the uncertainty bounds in \eqref{eq:load} became tighter and progressively more candidate vectors were rendered admissible. Even for loose voltage regulation limits of $\pm 10\%$ and tight load uncertainty, more than $20\%$ of the candidates in $\mcS$ violated \eqref{eq:volt}. 

The reduction from $\mcS$ to $\mcS_r$ via \eqref{eq:farkas} can be generalized. For example, limits on line and transformer flows can be expressed as linear functions of power injections and appended to \eqref{eq:volt}. Moreover, correlations in load forecasts across buses, or power factor limitations applied on a per-bus basis, both can be directly captured as linear inequalities and appended to \eqref{eq:load}. Finally, if the library has been reduced significantly so that $L<T$, the operator could broaden the voltage interval $[\underline{u},\overline{u}]$ and/or tighten the load uncertainty range in \eqref{eq:volt} if grid probing is still needed to recover non-metered loads.

\subsection{Finding Probing Setpoints with Most Diverse States}\label{subsec:S3}
Given the reduced library $\mcS_r=\{\bs_\mcM^\ell\}_{\ell=1}^L$ of implementable and network-compliant candidates, the last step is to select the $T$ candidates yielding the most diverse states. Recall that the system state $\bv^\ell$ related to probing injection $\bs_\mcM^\ell$ depends also on the unknown loads $\bs_\mcO$. Moreover, the dependence on both $\bs_\mcM^\ell$ and $\bs_\mcO$ is non-linear and implicit. The approximate LDF model of \eqref{eq:LDF} can help us circumvent these technical challenges. 

The Euclidean distance between the system states induced by injections $\bs_\mcM^\ell,\bs_\mcM^{\ell'}\in\mcS_r$ will be surrogated by the Euclidean distance between the approximate states of \eqref{eq:LDF} as
\begin{equation*}
\|\bv^\ell - \bv^{\ell'}\|_2\simeq \|\by_\ell - \by_{\ell'}\|_2
\end{equation*}
for all $\ell,\ell'=1,\ldots,L$. The latter simplifies as 
\begin{align*}
\|\by_\ell - \by_{\ell'}\|_2&=\left\|\begin{bmatrix}
\bK & \bL \\
\bM & \bN
\end{bmatrix}\left(
\begin{bmatrix}
\bs_\mcM^\ell \\
\bs_\mcO
\end{bmatrix}-
\begin{bmatrix}
\bs_\mcM^{\ell'} \\
\bs_\mcO
\end{bmatrix}
\right)\right\|_2\nonumber\\
&=\left\|\begin{bmatrix}
\bK\\
\bM
\end{bmatrix}
\left(\bs_\mcM^\ell-\bs_\mcM^{\ell'}\right)\right\|_2
\end{align*}
where we have exploited the linearity in \eqref{eq:LDF} together with the fact that non-metered loads remain roughly invariant during probing. We define the distance between $\bs_\mcM^\ell,\bs_\mcM^{\ell'}\in\mcS_r$ as
\begin{align}\label{eq:dist}
d(\ell,\ell')&:=\|\by_\ell - \by_{\ell'}\|_2^2\nonumber\\
&=(\bs_\mcM^\ell-\bs_\mcM^{\ell'})^\top
(\bK^\top\bK+\bM^\top\bM)
(\bs_\mcM^\ell-\bs_\mcM^{\ell'}).
\end{align}

Based on this metric, we would like to select a subset $\mcA$ of $T$ out of the $L$ candidate vectors in $\mcS_r$ so that the sum of their pairwise distances is maximized
\begin{align}\label{eq:MSD}
\max_{\mcA\subset\mcS_r}~&~\sum_{\ell\in\mcA}\sum_{\ell'\in\mcA} d(\ell,\ell')\\
\mathrm{s.to}~&~|\mcA|=T.\nonumber
\end{align}
The task in \eqref{eq:MSD} is known as the \emph{max-sum diversity} (MSD) problem, and appears frequently in information retrieval, computational geometry, and operations research~\cite{Cevallos16}. In fact, MSD can be reformulated as a binary quadratic program (QP) after introducing the $L\times L$ distance matrix $\bD$ with entries $D_{\ell,\ell'}:=d(\ell,\ell')$ as
\begin{subequations}\label{eq:MSD2}
\begin{align}
f^\star:=\max_{\bx\in\{0,1\}^L}~&~\bx^\top\bD\bx\label{eq:MSD2:cost}\\
\mathrm{s.to}~&~\bx^\top\bone=T.\label{eq:MSD2:c1}
\end{align}
\end{subequations}
Despite its simple form, the MSD task is NP-hard~\cite{Cevallos16}. However, thanks to the properties of $\bD$, the problem in \eqref{eq:MSD2} enjoys a polynomial-time approximate scheme (PTAS)~\cite{Cevallos16}.

Although $\bD$ is indefinite, the objective in \eqref{eq:MSD2:cost} can be shown to be concave under constraint \eqref{eq:MSD2:c1}. To see this, define the $2N\times L$ matrix $\tbY:=\left[\by_1~\cdots~\by_L\right]$ and use the definition of $d(\ell,\ell')$ to rewrite the objective of \eqref{eq:MSD2} as
\begin{align*}
f(\bx)&:=\bx^\top\bD\bx=\sum_{\ell=1}^L\sum_{\ell'=1}^L x_\ell x_{\ell'} D_{\ell,\ell'}\\
&=\sum_{\ell=1}^L\sum_{\ell'=1}^L x_\ell x_{\ell'} \left(\|\by_\ell\|_2^2+\|\by_{\ell'}\|_2^2-2\by_\ell^\top\by_{\ell'}\right)\\
&=\sum_{\ell'=1}^L  x_{\ell'} \|\by_{\ell'}\|_2^2\left(\sum_{\ell=1}^L x_\ell \right) + \sum_{\ell=1}^L  x_{\ell} \|\by_{\ell}\|_2^2\left(\sum_{\ell'=1}^L x_{\ell'} \right)\\
&-2\sum_{\ell=1}^L\sum_{\ell'=1}^L x_\ell x_{\ell'}\by_\ell^\top\by_{\ell'}\\
&=2T\bc^\top\bx -2\bx^\top \tbY^\top\tbY \bx.
\end{align*}
where $\bc:=\left[\|\by_1\|_2^2~\cdots~\|\by_L\|_2^2\right]^\top$. Since $\tbY^\top\tbY\succeq \bzero$, the objective $f(\bx)$ equals a concave quadratic function.

For moderate $L$ (a few hundreds), the task in \eqref{eq:MSD2} can be handled by a mixed-integer QP solver. For $T=2$, the MSD solution can be found by an exhaustive search. For larger $T$, we will use a randomized rounding approach, as adopted from \cite{Rag1987} in \cite[Remark~2]{Cevallos16}. The approach is briefly reviewed here for completeness. Its first step solves the relaxed problem
\begin{subequations}\label{eq:MSDr}
\begin{align}
\hbx:=\arg\min_{\bzero\leq \bx\leq \bone}~&~2\bx^\top \tbY^\top\tbY \bx- 2T\bc^\top\bx\label{eq:MSDr:cost}\\
\mathrm{s.to}~&~\bx^\top\bone=T.\label{eq:MSDr:c1}
\end{align}
\end{subequations}
Since the binary constraints of \eqref{eq:MSD2} are related to box constraints in \eqref{eq:MSDr}, it holds that $f(\hbx)\geq f^\star$. To construct a point $\tbx$ that is feasible for \eqref{eq:MSD2}, draw $L$-dimensional vectors $\{\tbx_i\}$ whose entries are independent Bernoulli random variables with mean $(1-\beta)\hbx$ for some $\beta>0$, say $\beta=0.1$. The so constructed binary vectors $\tbx_i$'s satisfy $\mathbb{E}[\tbx_i^\top\bone]=(1-\beta)T$ and $\mathbb{E}[\tbx_i^\top\bD\tbx_i]=(1-\beta)^2\hbx^\top\bD\hbx$. The purpose of scaling $\hbx$ by $(1-\beta)$ is to ensure $\tbx_i$'s are both feasible for \eqref{eq:MSD2} and yield relatively high cost with significant probability~\cite{Cevallos16}. 

Let us now comment on the complexity for designing probing setpoints. The first step described in Section~\ref{subsec:S1} is computationally inexpensive. The second step of Section~\ref{subsec:S2} involves solving the linear program in \eqref{eq:farkas} $K$ times, once for each candidate setpoint vector. The third step of Section~\ref{subsec:S3} entails solving the linearly-constrained quadratic program of \eqref{eq:MSDr}, whose complexity is cubic in the number of variables $L$. As detailed later in Section~\ref{sec:sim}, running this design process for the IEEE 34-bus feeder and $K=100$ candidate setpoints took $1-1.5$ min depending on $(M,O)$. The tests were run on a laptop computer using generic off-the-shelf solvers.

\begin{figure}[t]
	\centering
	\includegraphics[scale=0.4]{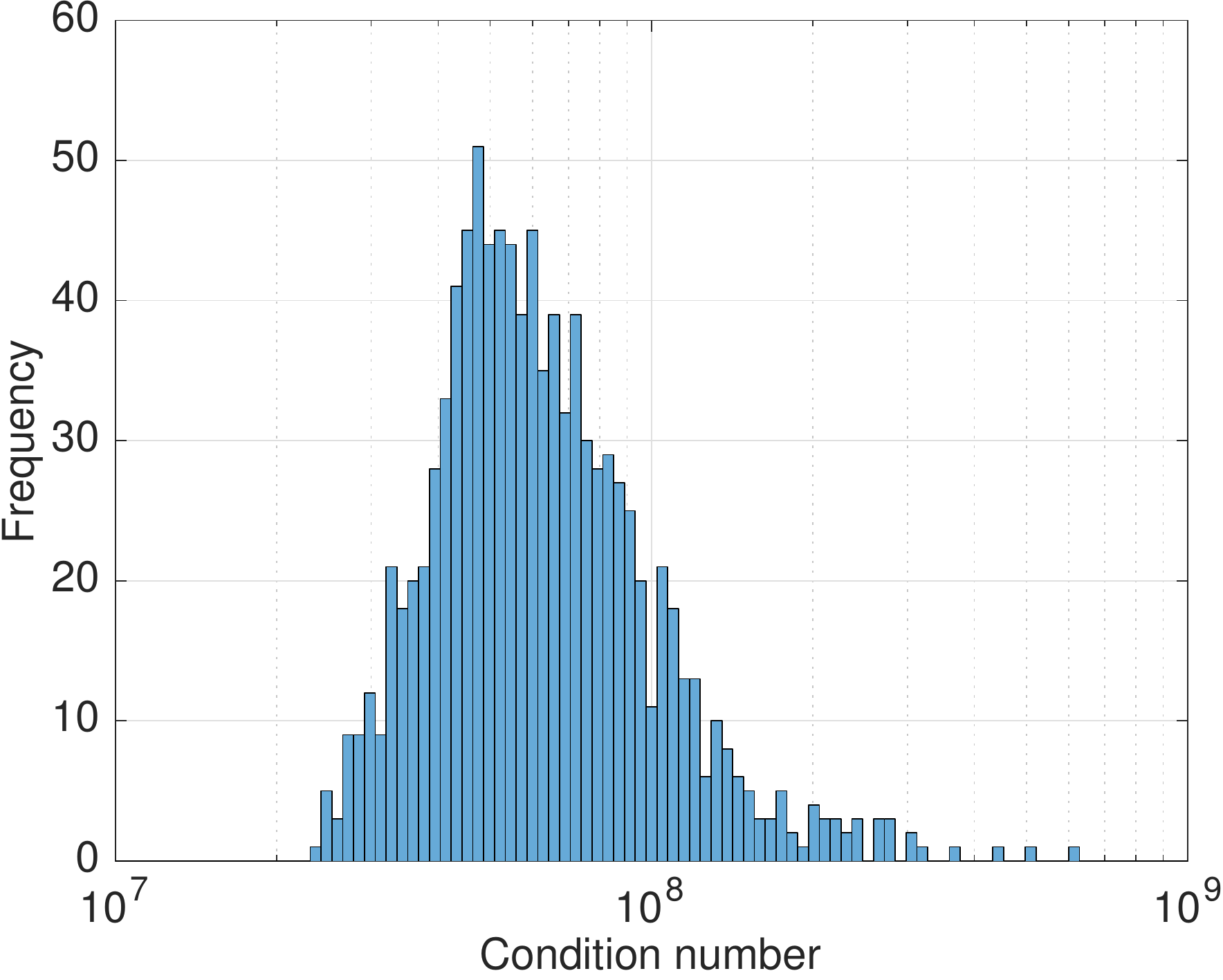}
	\caption{Histogram of condition numbers for the Jacobian matrix $\mathbf{J}(\{\mathbf{v}_t\}_{t=1}^4)$ obtained by randomly sampling quadruplets of $\bs_\mcM$'s from $\mcS_r$.}
	\label{fig:barplot}
\end{figure}

\begin{figure}[t]
	\centering
	\includegraphics[scale=0.32]{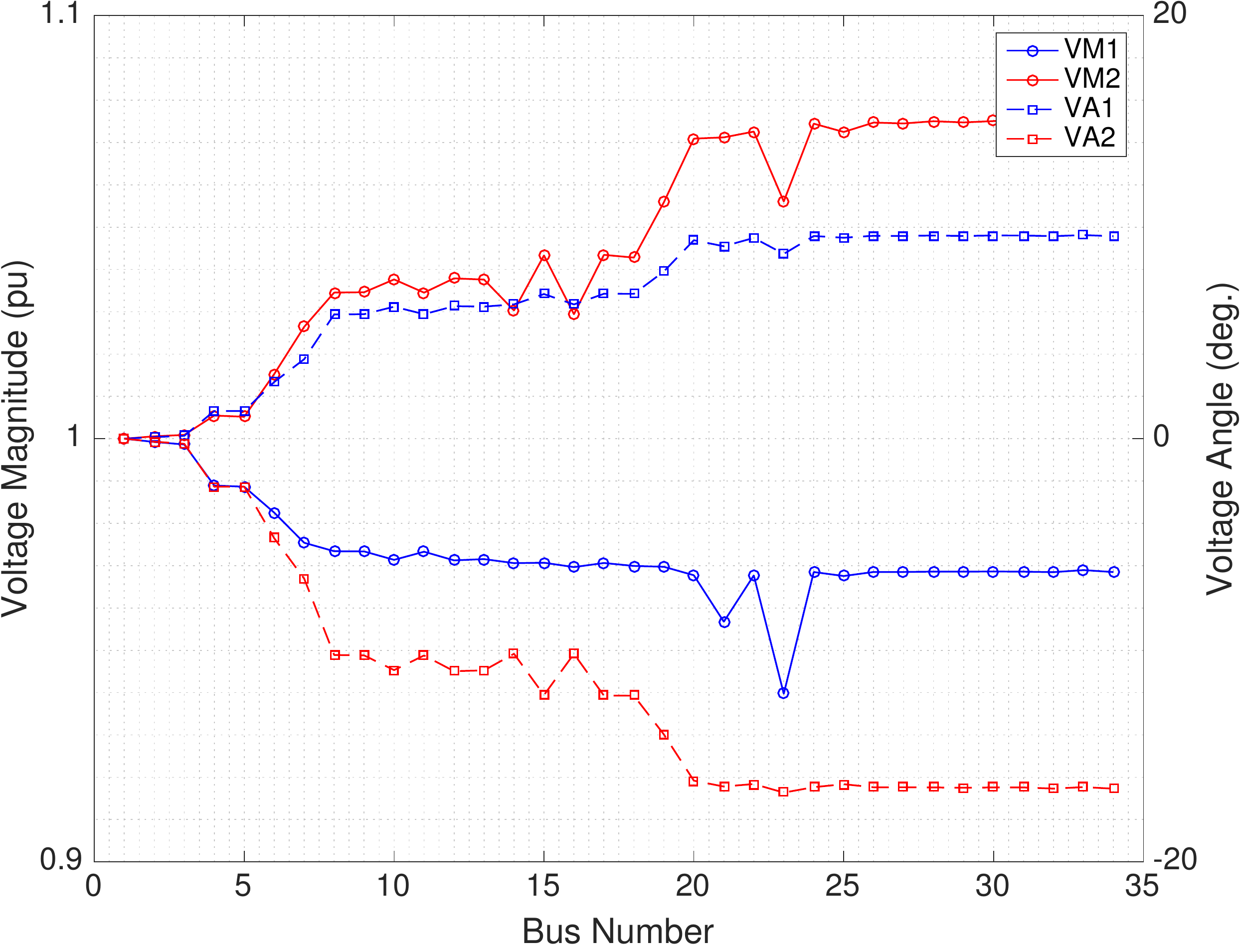}
	\caption{Probing setpoint design on the IEEE 34-bus grid for $T=2$ and $O=6$. The blue (red) lines correspond to the system states induced by the first (second) probing setpoints $\{\bs_\mcM^t\}_{t=1}^2$. Solid lines depict voltage magnitudes, and dashed lines voltage angles.}
	\label{fig:clustT2}
\end{figure}

To justify the need for this third step in probing design, we conducted a test on the IEEE 34-bus feeder for $T=4$ and $O=6$. For this test, load uncertainty was confined within a factor of $\pm 1$ times the nominal loads. Given library $\mcS$ of randomized injections drawn from $\pm 0.2$~pu and obeying \eqref{eq:appc1}--\eqref{eq:appc3}, we constructed the reduced library $\mcS_r$ based on \eqref{eq:farkas} for $[\underline{u},\overline{u}]=[0.90,1.10]$. We then solved \eqref{eq:MSDr} and followed the randomized rounding process to construct 100 binary $\tbx_i$'s. We evaluated the cost $f(\tbx_i)$ for those $\tbx_i$'s satisfying $\tbx_i^\top\bone=T$, and returned the $\tbx_i$ yielding the largest cost. The condition number of the Jacobian matrix evaluated at the so obtained $\tbx_i$ was $3\cdot 10^6$. We also calculated the condition number of the Jacobian matrix evaluated at random candidate quadruplets in $\mcS_r$. The latter condition numbers ranged within $10^7-10^9$; see Fig.~\ref{fig:barplot}. Hence, albeit MSD adds computational complexity, it is an important part of the probing design process. 

To show that this MSD step provides diverse system states, we performed another test on the IEEE 34-bus feeder under the same setup, but for $T=2$. The states induced by the designed probing and nominal loads are shown in Figure~\ref{fig:clustT2}. 

Upon solving \eqref{eq:MSD} near optimally, we have obtained $T$ probing injection vectors $\{\bs_\mcM^t\}_{t=1}^T$ that are: \emph{i)} implementable by inverters; \emph{ii)} network-compliant; and \emph{iii)} yield diverse system states. In the process of probing design, the first step (Section~\ref{subsec:S1}) operates on the entries of $\bs_\mcM^t$'s; the second step (Section~\ref{subsec:S2}) considers each vector $\bs_\mcM^t$ as a whole; and the third step (Section~\ref{subsec:S3}) accounts for the joint effect of probing injections $\{\bs_\mcM^t\}_{t=1}^T$.

\section{Solving the P2L tasks}\label{sec:solvers}
Recall that the P2L task with phasor data involves solving the set of non-linear equations
\begin{subequations}\label{eq:PF12}
	\begin{align}
	u_n(\bv_t)&=u_n^t &\forall n\in\mcM, t \in \mcT \label{eq:PF12u}\\
	\theta_n(\bv_t)&=\theta_n^t &\forall n\in\mcM, t \in \mcT\label{eq:PF12t}\\
	p_{n}(\bv_t)&=p_{n}^t &\forall n\in\mcM, t \in \mcT \label{eq:PF12p}\\
	q_{n}(\bv_t)&=q_{n}^t &\forall n\in\mcM, t \in \mcT \label{eq:PF12q}\\
    p_{n}(\bv_t)&=p_{n}(\bv_{t+1}) &\forall n\in\mcO,t \in \mcT' \label{eq:PFC2p}\\
	q_{n}(\bv_t)&=q_{n}(\bv_{t+1}) &\forall n\in\mcO,t \in \mcT' \label{eq:PFC2q}
	\end{align}
\end{subequations} 
where $\bv_t$'s are the system states across $\mcT=\{1,\ldots,T\}$; $\{(u_n^t,\theta_n^t,p_n^t,q_n^t)\}_{n\in\mcM}$ are the probing data collected at time $t$; \eqref{eq:PF12u}--\eqref{eq:PF12q} are the $4MT$ metering equations; and \eqref{eq:PFC2p}--\eqref{eq:PFC2q} are the $2O(T-1)$ coupling equations with $\mcT':=\{1,\ldots,T-1\}$. For the P2L task with non-phasor data, the angle information in \eqref{eq:PF12t} is unavailable. 

Having characterized the local identifiability for the P2L tasks in Part I, this section presents solvers for tackling P2L. If grid specifications are noiseless and the power injections in $\mathcal{O}$ remain unaltered during probing, the P2L tasks boil down to solving the equations in \eqref{eq:PF12}. The latter can be tackled by adopting the semidefinite program (SDP)-based solvers developed in \cite{ZhuGia12}, \cite{KlauberZhu15}, \cite{MLB15}, \cite{MALB16}. Here we will skip the details, which can be found in \cite{BKV17} for $T=2$, and outline the P2L solver for noisy data.

Probing data are inexact due to measurement noise and modeling inaccuracies in the metering equations of \eqref{eq:PF12u}--\eqref{eq:PF12q}. To account for small fluctuations in non-metered loads during probing, a noise term is added to the RHS of the coupling equations in \eqref{eq:PFC2p}--\eqref{eq:PFC2q}. To cope with noisy data, we extend the penalized SDP-based state estimator of~\cite{MALB16} to the P2L setting as follows
\begin{subequations}\label{eq:CPSSE}
	\begin{align}
	\min~ &\alpha \sum_{t=1}^{T}{\trace(\mathbf{M}\bV_t)} + \sum_{t=1}^{T}\sum_{k=1}^{3M} f_{k}({\epsilon_{k}^t}) + \sum_{t=1}^{T-1}\sum_{l=1}^{2O} f_{l}({\xi_{l}^t}) \label{eq:CPSSE:cost}\\
\mathrm{over}~&~		\bV_t \succeq \bzero, \{\epsilon_k^t\}_{k=1}^{3M},\quad t\in\mcT\\
~&~\{\xi_l^t\}_{l=1}^{2O},\quad t\in\mcT'\\
	\mathrm{s.to} ~&~\trace(\mathbf{M}_k\bV_t)+{\epsilon_{k}^t}=\hat{s}^t_k, \quad k=1:3M, t\in\mathcal{T} \label{eq:CPSSE:1}\\	~&~\trace(\mathbf{M}_l\bV_t)=\trace(\mathbf{M}_l\bV_{t+1})+{\xi_{l}^t}, \quad l=1:2O, t\in\mathcal{T}' \label{eq:CPSSE:O}
	\end{align}
\end{subequations}
where the given matrices $\mathbf{M}_k$ depend on $\mathbf{Y}$; see \cite{KlauberZhu15}, \cite{MLB15}. 

The matrix variables $\bV_t\succeq\bzero$ have been obtained upon relaxing the rank-one constraint $\bV_t=\tilde{\bv}_t\tilde{\bv}_t^H$ on the original system states for $t\in\mcT$. The measurements $\hat{s}_k^t$ relate to state $\mathbf{v}_t$ in \eqref{eq:CPSSE:1}; and the constraints in \eqref{eq:CPSSE:O} couple the $T$ states. The auxiliary variables $\epsilon_k^t$ can be substituted from \eqref{eq:CPSSE:1}--\eqref{eq:CPSSE:O} into the objective of \eqref{eq:CPSSE}; they are introduced here only to simplify notation. The data fitting penalties $f_k$ can be either a weighted squared or absolute value, that is
\begin{subequations}
	\begin{align*}
	f_k(\epsilon_k^t) &=\left(\frac{\epsilon_k^t}{\sigma_k^2}\right)^2= \frac{\left(\hat{s}^t_k-\trace(\mathbf{M}_k\bV_t)\right)^2}{\sigma_k^2} \quad \text{or} \\ 
	f_k(\epsilon_k^t) &=\frac{|\epsilon_k^t|}{\sigma_k}=\frac{|\hat{s}^t_k-\trace(\mathbf{M}_k\bV_t)|}{\sigma_k}
	\end{align*}
\end{subequations}
with different $\sigma_k$'s depending on the uncertainty of the $k$-th datum. Likewise, the auxiliary variables $\xi_l^t$'s capture variations of non-metered loads and are penalized through $f_l$'s, which are defined as $f_k$'s.

The first summand in \eqref{eq:CPSSE:cost} corresponds to a regularizer promoting rank-one minimizers for $\bV_t$; a practical choice sets $\bM=\bG$ as suggested in \cite{MLB15}. The second and third summands in \eqref{eq:CPSSE:cost} are data-fitting terms. The tuning parameter $\alpha>0$ governs the balance between the regularizer and the data-fitting terms: For $\alpha=0$, the P2L cost involves only the data-fitting terms; whereas for increasing $\alpha$, more emphasis is placed on the regularizer~\cite{MALB16}. If one or more of the minimizers $\bV_t^\star$ of \eqref{eq:CPSSE} is not rank-one, the heuristic for constructing a system state $\bv_t^\star$ proposed in~\cite{MALB16} is used. 

Additional constraints can be added to strengthen the SDP relaxation. For example, if non-metered buses are known to host exclusively loads, the constraints  $\trace(\mathbf{M}_l\bV_t)\leq0$ for $l=1:2O$, and $t \in \mathcal{T}$ can be appended to \eqref{eq:CPSSE}. Additional information on loads, such as the uncertainty range of \eqref{eq:load}, can be readily incorporated. As in \cite{MALB16}, if bus $n$ is known to be a zero-injection bus, then $\tilde{i}_n=\be_n^\top\bY\tbv$ has to be zero. Therefore, the constraint $\tbv\tilde{i}_n^\star=\bV \bY^\star\mathbf{e}_n=\bzero$ can be added. 

Given phasor data, the metering equations corresponding to voltage magnitudes can be dropped. If the vectors of voltage phasors $\{\tbv_t\}_{t=1}^T$ are included as optimization variables, the direct measurements on the voltage phasors of $\mcM$ can be simply expressed as
\begin{equation}\label{eq:PMU1}
\tilde{v}_{t,k}+{\epsilon_{k}^t}=\hat{s}^t_k, \quad k=1:M,~t\in\mcT.
\end{equation}
To capture the dependence between $\tbv_t$ and $\tbV_t$, the non-convex constraint
\[\rank\left(	\begin{bmatrix} 
	\bV_t & \tilde{\bv}_t \label{eq:PMUb} \\
	\tilde{\bv}_t^H & 1
	\end{bmatrix}\right)=1\] 
can be surrogated by the next SDP constraint as in \cite{ZhuGia12}
\begin{equation} \label{eq:PMU2}
	\begin{bmatrix} 
	\bV_t & \tilde{\bv}_t\\
	\tilde{\bv}_t^H & 1
	\end{bmatrix}
	\succeq \bzero, \quad t\in\mcT.
\end{equation}
Since the $\tbv_t$'s are optimization variables now, there is no need to use the heuristic of \cite{MALB16} to recover the system states.

\section{Numerical Tests}\label{sec:sim}
The topological observability criteria for the P2L task and the SDP-based solvers were numerically tested using the IEEE 34-bus feeder. The original multi-phase grid was converted to an equivalent single-phase grid~\cite{GLTL12}. The numerical tests were run on a 2.7 GHz Intel Core i5 laptop computer with 8 GB RAM using the Sedumi solver on YALMIP and MATLAB~\cite{sedumi}, \cite{YALMIP}.

\subsection{Numerical Observability}\label{subsec:Jacobian}
Since Theorems~1 and 2 of Part I rely on the sparsity pattern rather than the exact values of $\bJ\left(\{\bv_t\}\right)$, we evaluated $\bJ\left(\{\bv_t\}\right)$ for 1,000 random state sequences $\{\bv_t\}_{t=1}^T$. The scenarios of phasor and non-phasor data were tested under four probing setups. For each setup, the placement of non-metered $\mcO$ and probing buses $\mcM$ were fixed. We generated 1,000 random state sequences by randomly drawing voltage magnitudes in the range $[0.90,1.10]$ per unit and voltage angles in the range $[-1.5, 1.5]\degree$. Assuming non-phasor data first, the following four setups were constructed according to the condition of Th.~2 of Part I:
\begin{itemize}
\item Setup A meets the condition for $O=16$ and $T=2$.
\item Setup B meets the condition for $O=6$ and $T=2$.
\item Setup C does not meet the condition for $O=16$ and $T=2$, but it does for $T=4$.
\item Setup D does not meet the condition for $O=6$ and $T=2$, but it does for $T=4$.
\end{itemize}
The same setups were considered for phasor data. As discussed in Part I, setups A and B meet also the condition of Theorem~1. Additionally, setups C and D were constructed such that they meet the condition of Theorem~1 for $T=2$.

\begin{figure}[t]
	\centering
	\includegraphics[scale=0.4]{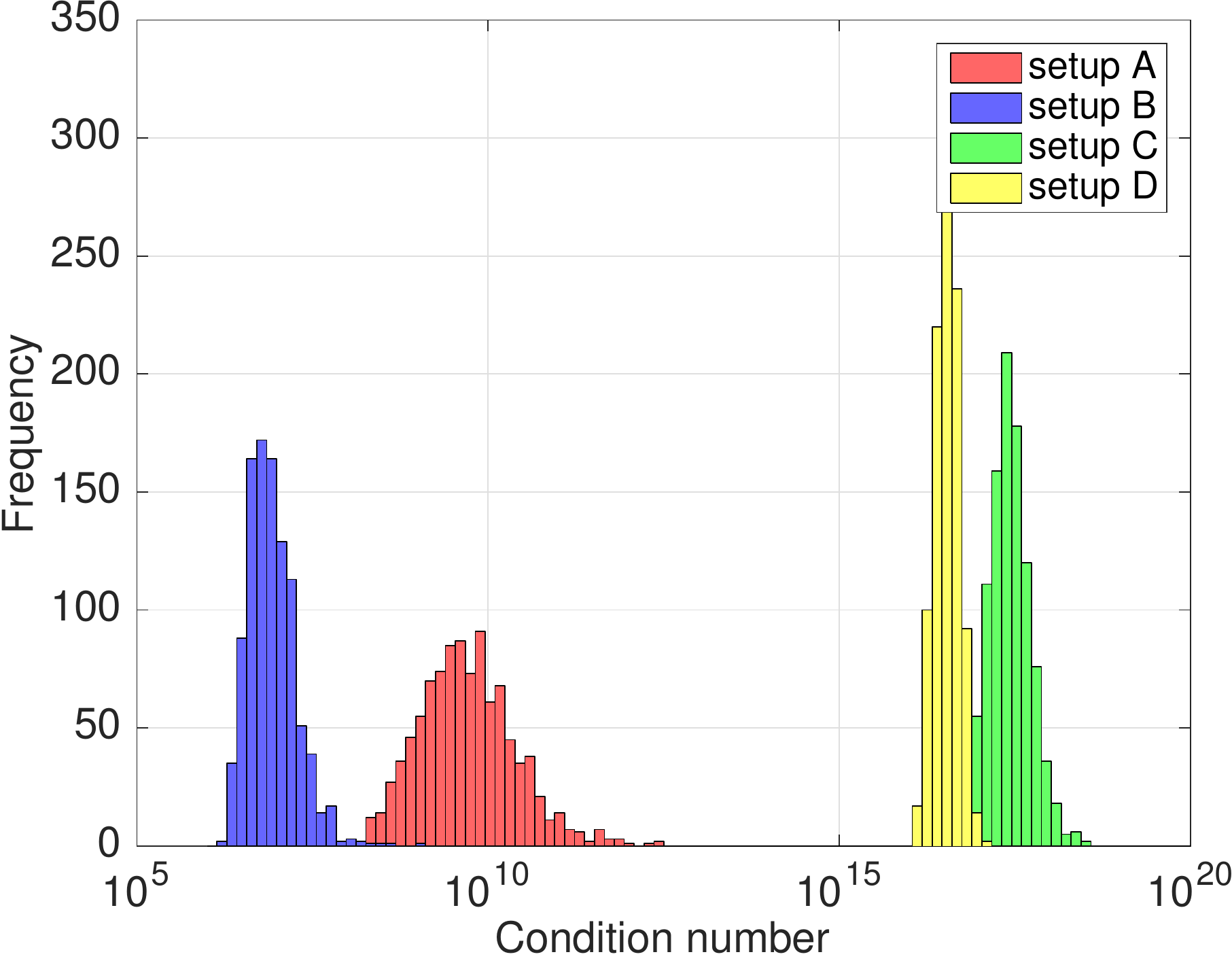}\\ \vspace*{1em}
	\includegraphics[scale=0.4]{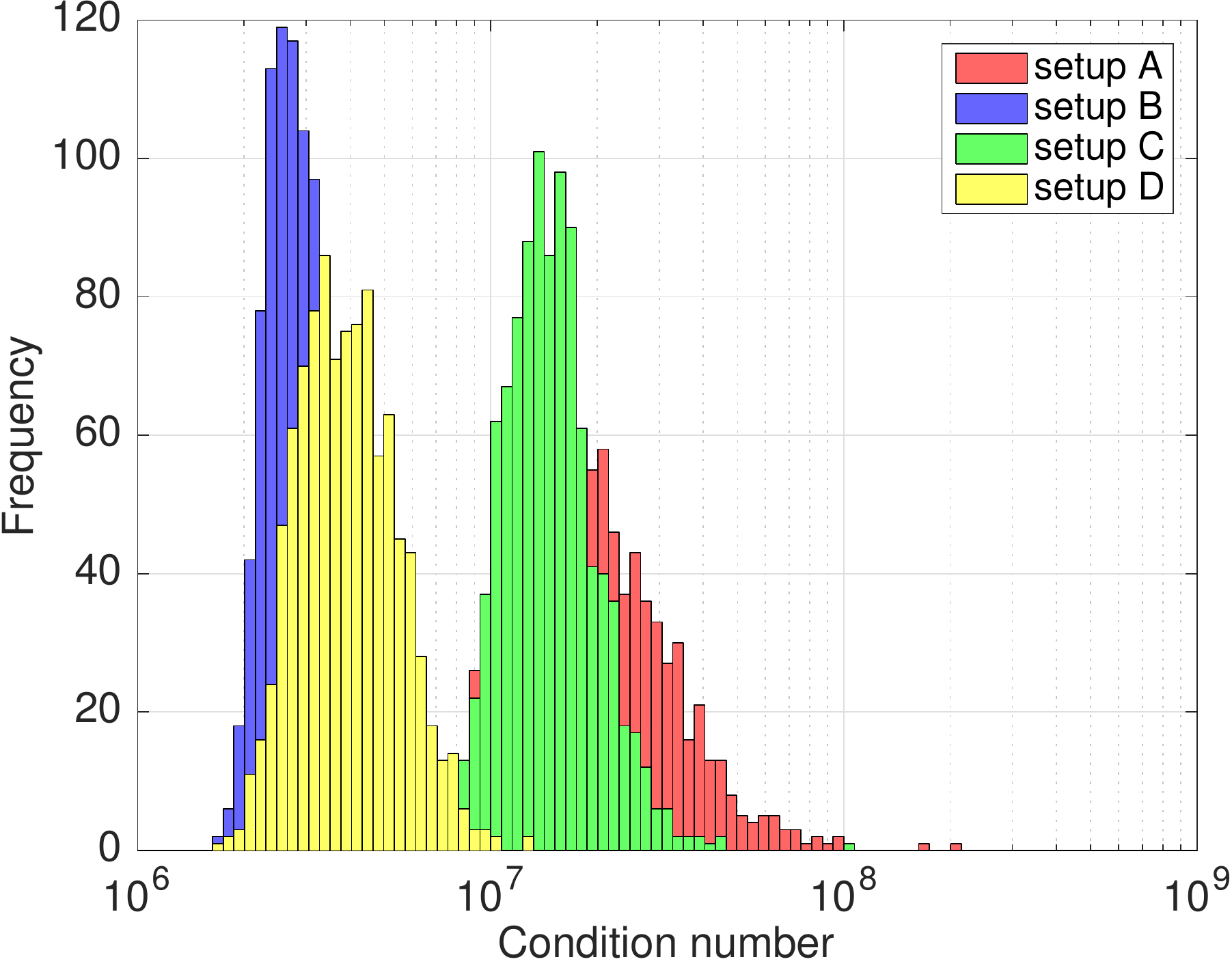}
	\caption{Histograms of the condition numbers for the P2L Jacobian matrices with non-phasor data for $T=2$ (top) and $T=4$ (bottom) probing actions.}
	\label{fig:Jacob}
\end{figure}

\emph{Non-phasor data:} Figure~\ref{fig:Jacob} depicts the condition number histograms obtained under the four setups for $T=2$ and $4$. Under setups A and B, although the dimensions of $\bJ\left(\{\bv_t\}\right)$ increase with $T$, the condition numbers did not. In fact, the condition number was sometimes reduced, especially in networks with large $\mcO$. For setups C and D, there was a significant shift in the histograms from $T=2$ to $T=4$, which validates Theorem~2. By and large, the condition number improves for decreasing $O$ and increasing $T$. Hence, when more loads are to be recovered, longer probing periods should be used. Of course, longer probing periods may violate the stationarity assumption on loads.

\begin{figure}[t]
	\centering
	\includegraphics[scale=0.4]{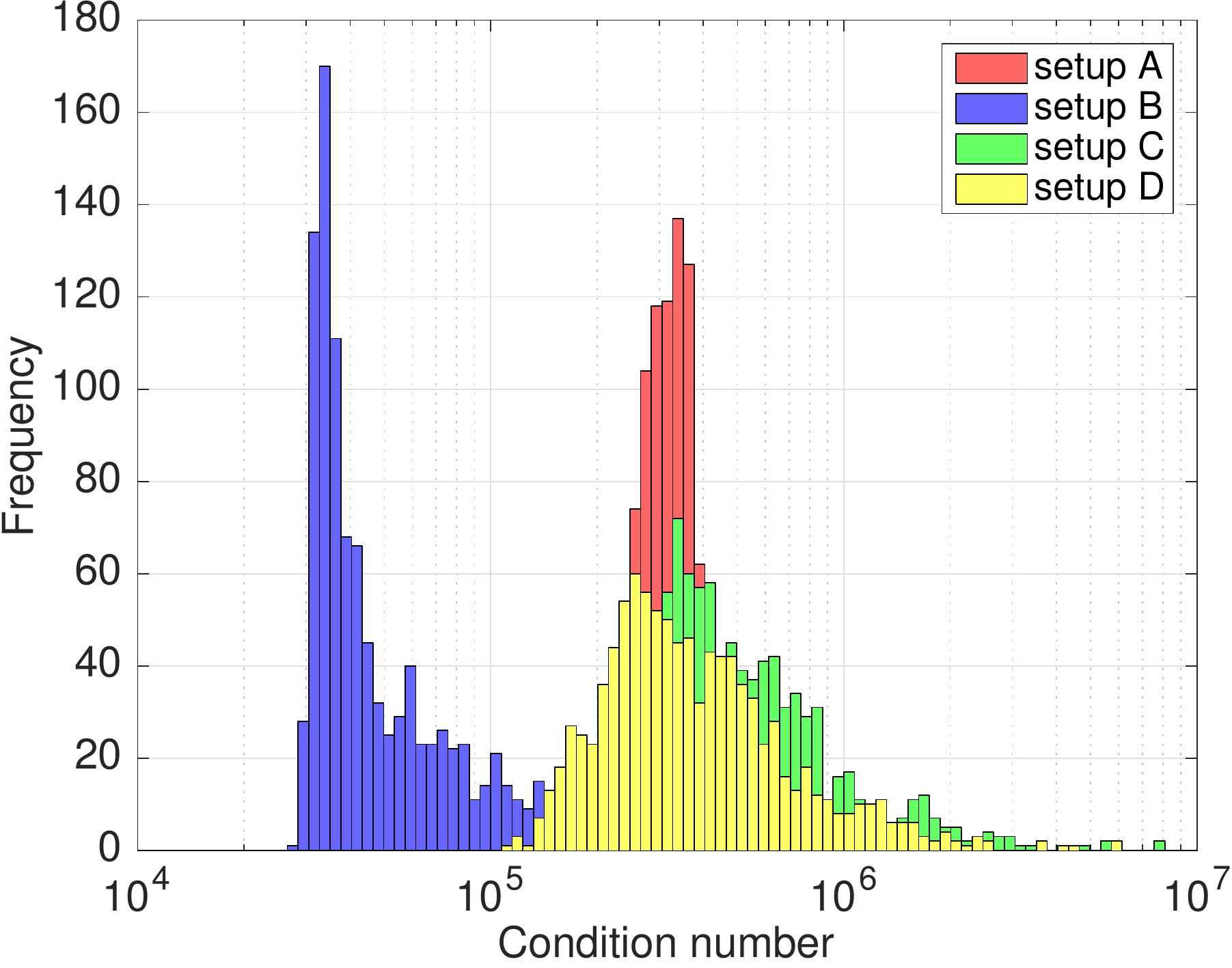}\\ \vspace*{1em}
 ~\includegraphics[scale=0.4]{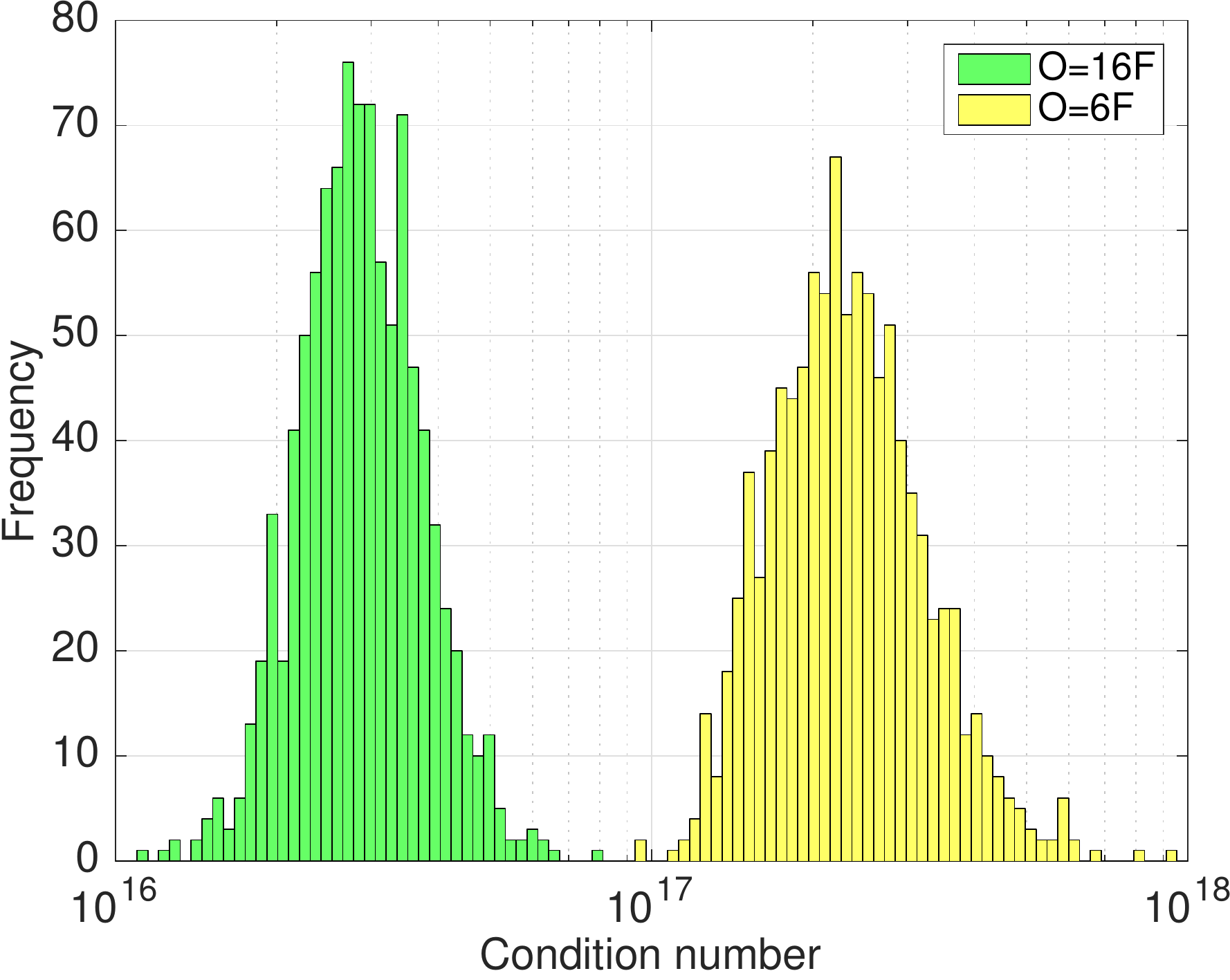}
	\caption{Histograms of the condition numbers for the P2L Jacobian matrices with phasor data for $T=2$ (top) and $T=4$ (bottom) probing actions.}
	\label{fig:Jacobpmu}
\end{figure}

\emph{Phasor data:} Figure~\ref{fig:Jacobpmu} displays the condition number histograms of $\bJ\left(\{\bv_t\}\right)$ again for $T=2$ and $4$. As expected, due to the value added of phasor data, the condition numbers decrease significantly. In addition, setups C and D that failed for $T=2$ with non-phasor data, become successful with phasor probing data. The tests corroborate the criteria of Th.~1. The bottom panel of Figure~\ref{fig:Jacobpmu} displays the condition number histograms under the following two setups that did not satisfy the condition of Th.~1: \emph{i)} for $T=4$ and $O=6)$ (yellow histogram); and \emph{ii)} for $T=4$ and $O=16$ (green histogram).  

\begin{figure}[t]
	\centering
	\includegraphics[scale=0.4]{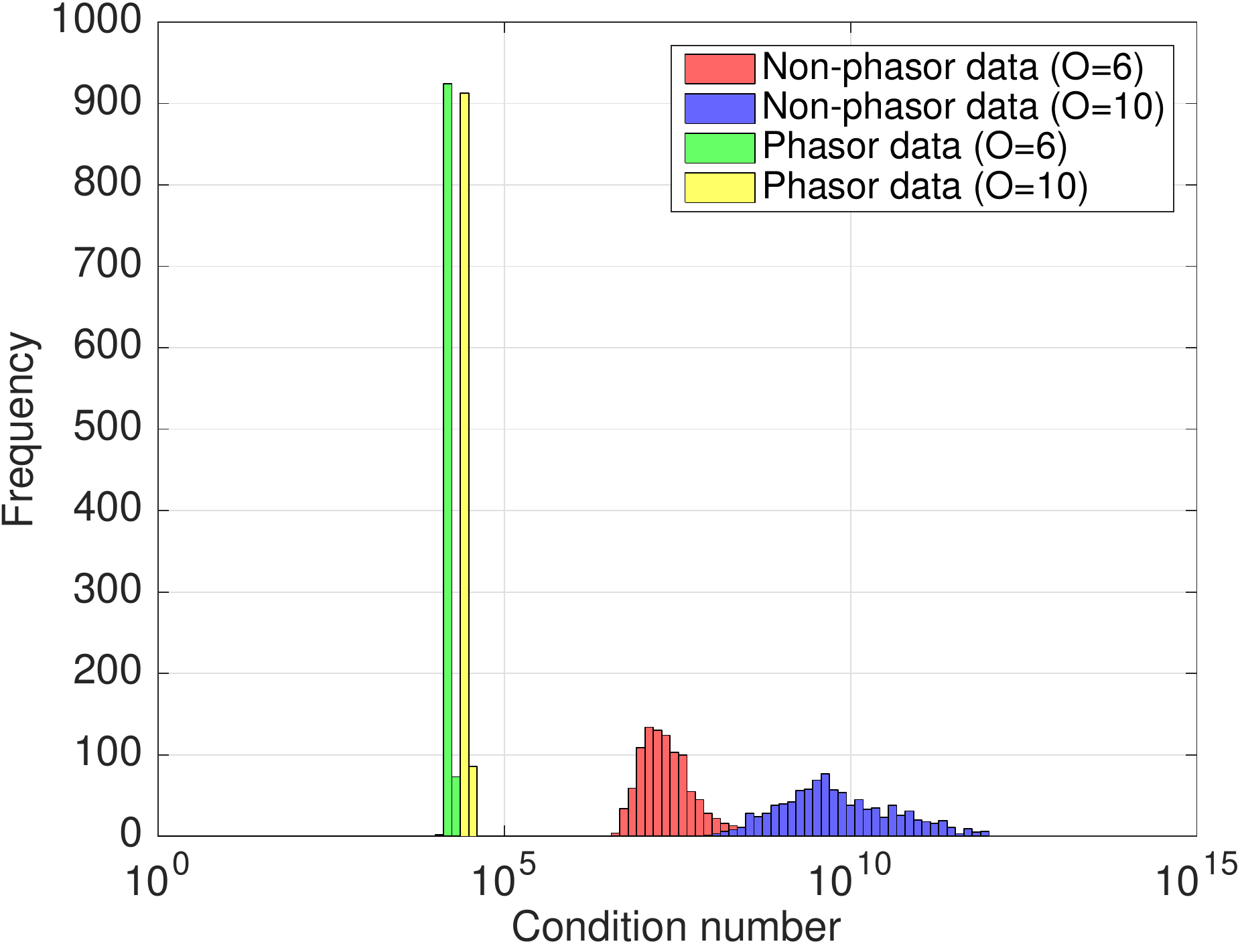}
	\caption{Histograms of the condition numbers for the Jacobian matrix $\bJ\left(\bv_1\right)$ with phasor and non-phasor data for single-slot probing $(T=1)$.}
	\label{fig:T1P}
\end{figure}

\emph{Single-slot probing scenario:} We also tested the special case of $T=1$, where one fixes voltages and injections on a subset of buses $\mcM$ and tries to find the loads at the remaining buses $\mcO$. This setup is pertinent to learning ZIP loads as discussed in Part I. We tested two fixed placements of non-metered buses that met the conditions of Theorems~3 and 4, respectively. We then evaluated $\bJ\left(\bv_1\right)$ at $1,000$ random system states. Figure~\ref{fig:T1P} shows the histograms for the condition numbers of $\bJ\left(\bv_1\right)$. Bus placements that did not meet the criteria of Th.~3 and 4 exhibited condition numbers similar to those at the bottom panel of Figure~\ref{fig:Jacobpmu}.

The condition number of the Jacobian matrices in PSSE tasks for transmission systems is known to depend heavily on the specification set \cite{Clements83}, \cite{Baldick01}: A larger number of voltage magnitude and line flow measurements tends to yield a lower condition number. It is thus expected that adding line flow measurements would improve load and state estimation.

\subsection{SDP-based P2L} \label{sec:SDP}
Given noisy specifications, the P2L tasks were tackled using actual data and the SDP-based solver of \eqref{eq:CPSSE}--\eqref{eq:PMU2}. The loads on the IEEE 34-bus grid were taken from the Pecan Street dataset~\cite{pecandata}, between 10:00 a.m. and 01:40 p.m. on January 1, 2013, and in 10-minute intervals. Load sequences were scaled so that the peak active load over the tested period was $0.5$ pu. Lacking values for reactive loads, a lagging power factor of 0.9 was simulated for all loads. 

\begin{figure}[t]
	\centering
	\includegraphics[scale=0.4]{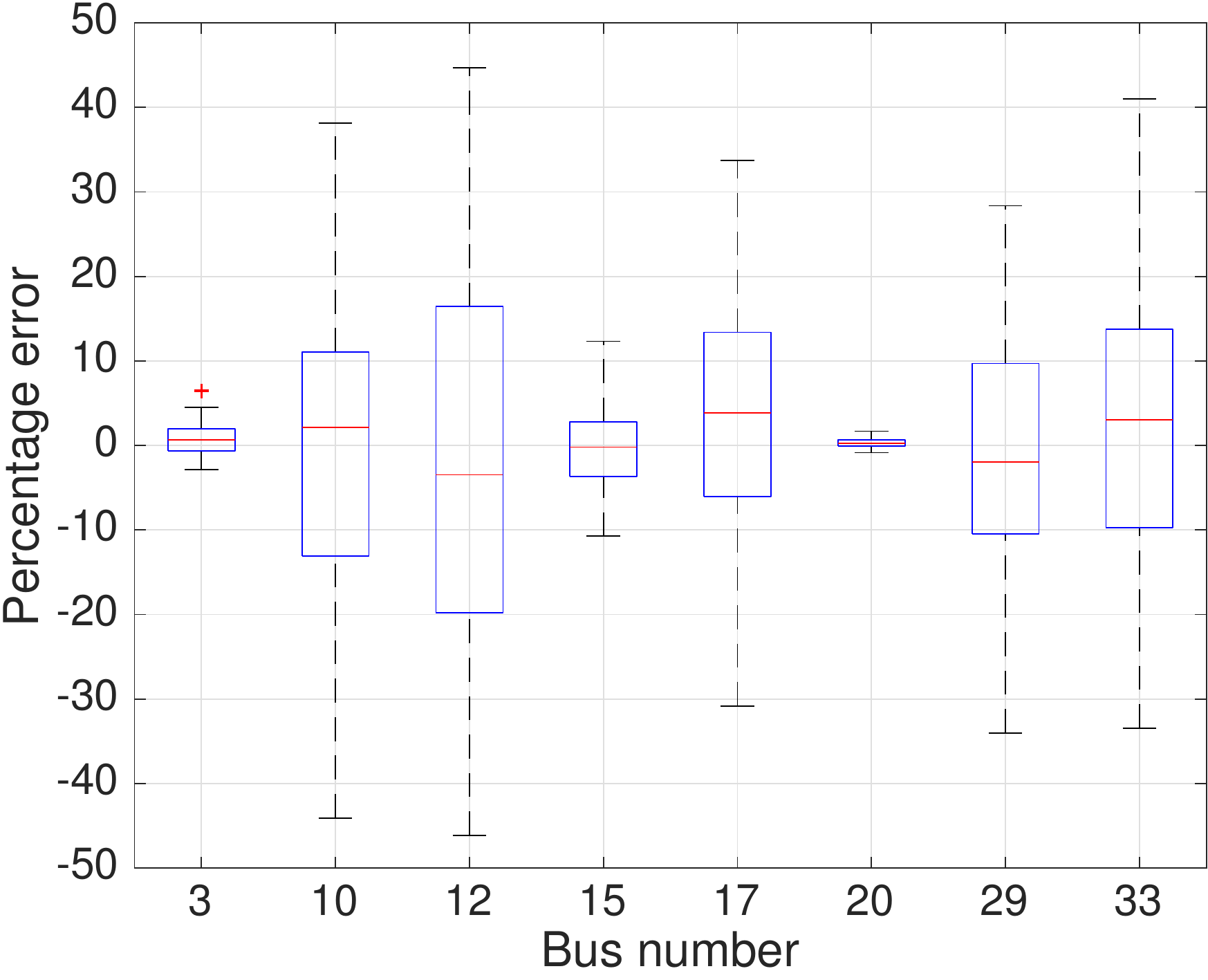}\vspace*{1em}
	\includegraphics[scale=0.4]{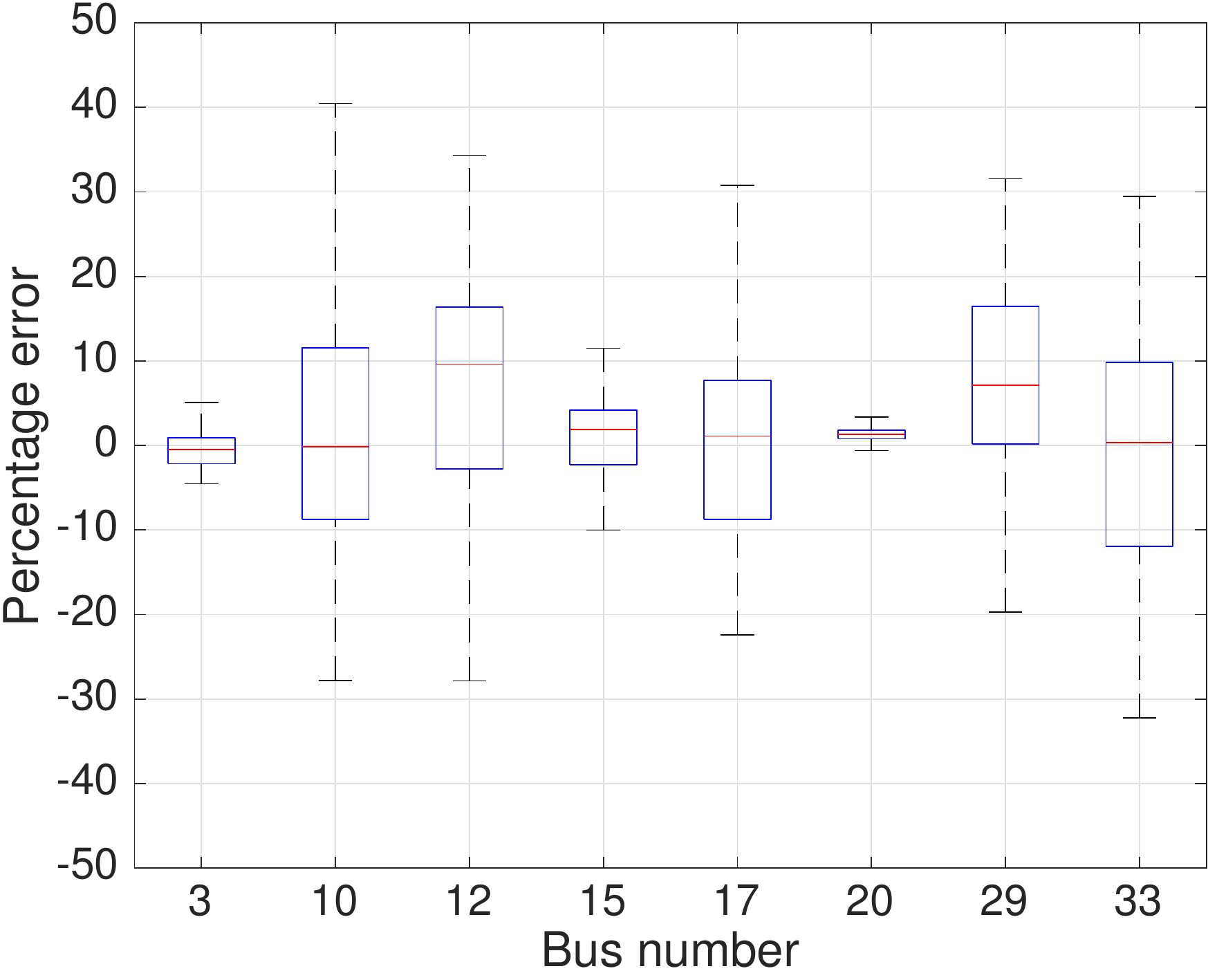} 
	\caption{Percentage error in active power injection estimates with phasor data for $T=4$ and $O=8$ without MSD (top) and with MSD (bottom).}
	\label{fig:MSD}
\end{figure}

To simulate probing injections at buses in $\mcM$, we first created a data library $\mcS$ of $K=100$ randomized injection vectors as described in Section~\ref{subsec:S1} for $\bar{p}_n=0.2$~pu. The library $\mcS$ was then reduced to $\mcS_r$ to ensure that voltage magnitudes lie within $[0.90,1.10]$~pu for non-metered loads within $[\mathbf{0}, 2\bs_{\mcO}]$ as described in Section~\ref{subsec:S2}. For all tests, the regularization parameter was set to $\alpha=20,000$, and the functions $f_k$ and $f_l$ in \eqref{eq:CPSSE} were selected as the WLS costs. To simulate measurement noise, the probing data recorded for an actual quantity $x$ (e.g., voltage magnitude or power injection) was modeled as $\hat{x}=x(1+\epsilon)$, where $\epsilon$ is a zero-mean Gaussian random variable. The variance $\sigma^2$ of $\epsilon$ was selected to yield the desired value of signal-to-noise ratio (SNR)
\begin{equation}\label{eq:SNR}
10\cdot \log_{10}\frac{x^2}{\mathbb{E}[x^2\epsilon^2]}=-20\cdot \log_{10}\sigma.
\end{equation}
This variance is the same variance appearing in \eqref{eq:CPSSE} as $\sigma_k^2$. Likewise, to capture small load variations, non-metered loads were simulated by perturbing their nominal value $p_n$ as $\hat{p}_n^t=(1+\epsilon_n)p_n$ for $t\in\mcT$, and similarly for $q_n$'s.

To check whether the MSD step of Section~\ref{subsec:S3} improves estimation, we tested P2L with and without this step. The test considered $100$ Monte Carlo realizations for the loads at 10:00 a.m. The P2L task was run for $T=4$, $O=8$, and using phasor data. The SNR values were set respectively to $80$ and $60$~dB for metered and non-metered buses. PMUs are expected to have such high accuracy~\cite{Frigo16}. The range of percentage errors was reduced from $[-50,+50]\%$ to $[-30, +40]\%$ by selecting the $T$ most diversifying setpoints. 

\begin{figure}[t]
	\centering
	\includegraphics[scale=0.4]{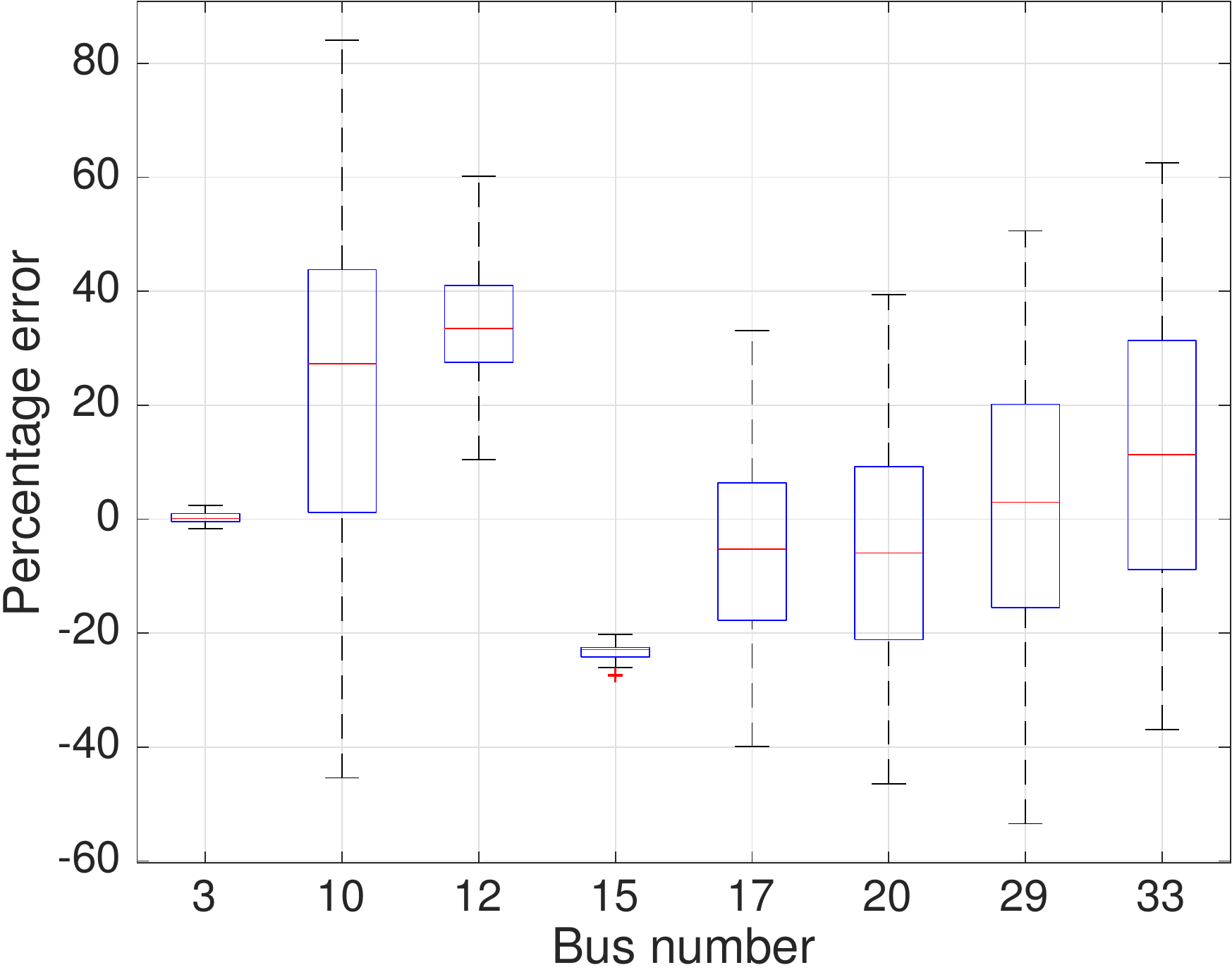} 
	\caption{Percentage error in active power injection estimates with non-phasor data for $T=4$ and $O=8$.}
	\label{fig:T6P}
\end{figure}

\begin{figure}[t]
	\centering
	\includegraphics[scale=0.4]{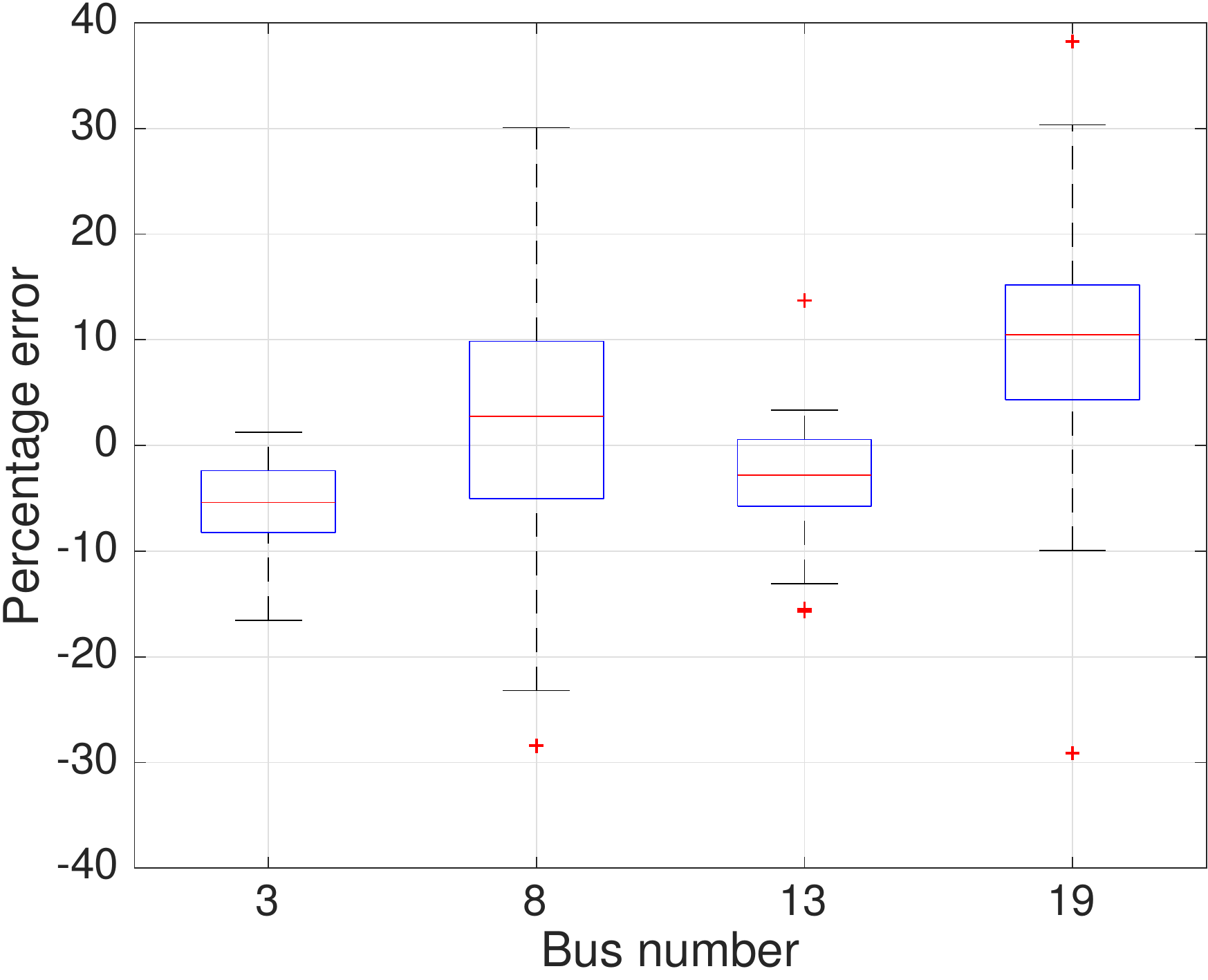}\\ \vspace*{1em}
	\includegraphics[scale=0.4]{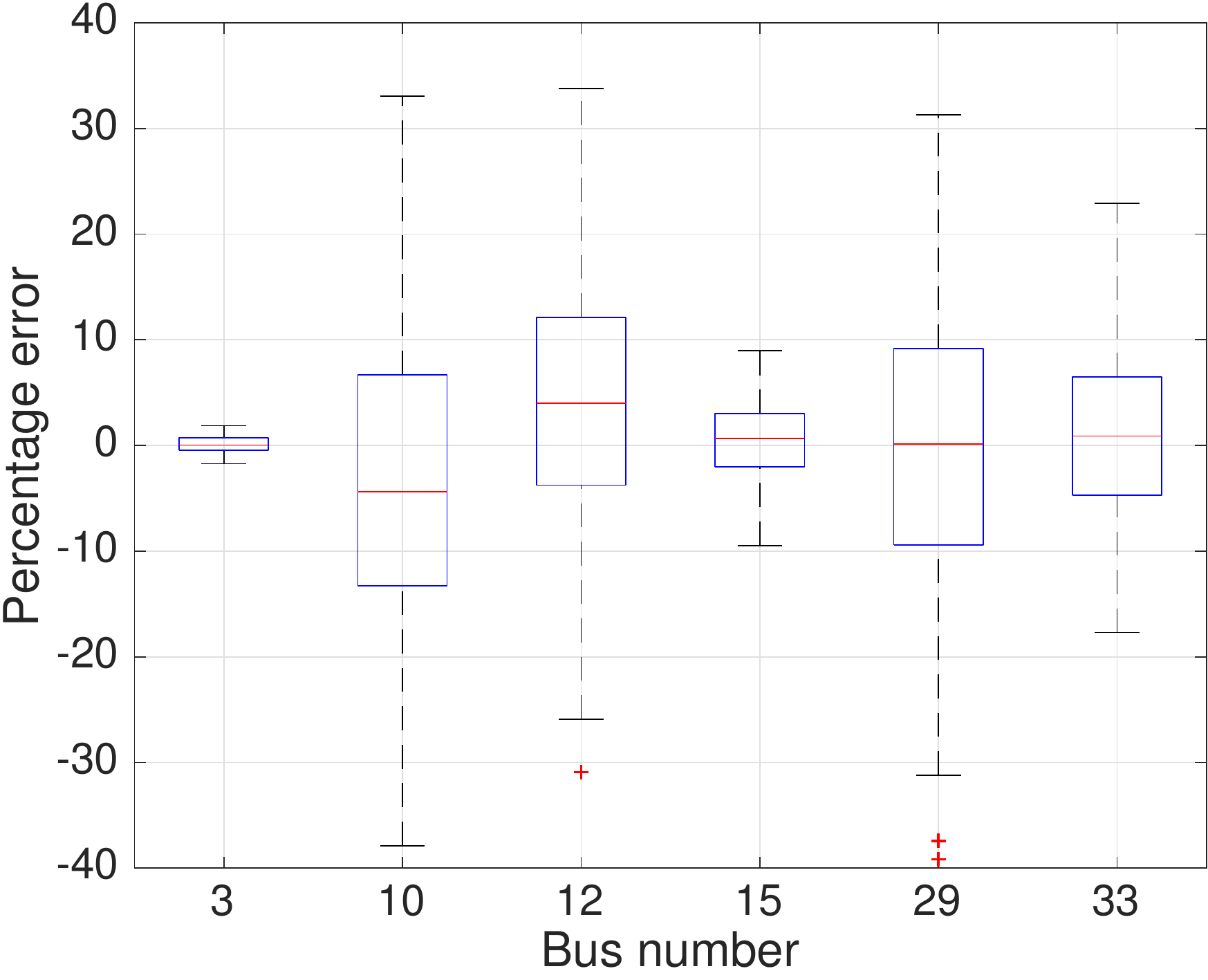}
	\caption{Percentage error in active power injection estimates for $T=1$ with non-phasor data and $O=4$ (top); and with phasor data and $O=6$ (bottom).}
	\label{fig:T1}
\end{figure}

To verify the improvement of using phasor over non-phasor probing data, we repeated the previous MSD setup but now for non-phasor data. The obtained percentage errors are depicted in Figure~\ref{fig:T6P} and are of worse accuracy compared to those in the bottom panel of Figure~\ref{fig:MSD}. We also tested the single-slot probing scenario of $T=1$ under slightly different probing setups for (non)-phasor data. Figure~\ref{fig:T1} illustrates the statistics of the obtained percentage errors.

\begin{remark}[]\label{re:}
Based on the numerical tests, we have observed that load estimates generally improve when: a) the MSD step is implemented; b) phasor data are utilized; c) the duration $T$ is increased; and d) $O$ is decreased.
\end{remark}

\begin{figure}[t]
	\centering
	\includegraphics[scale=0.39]{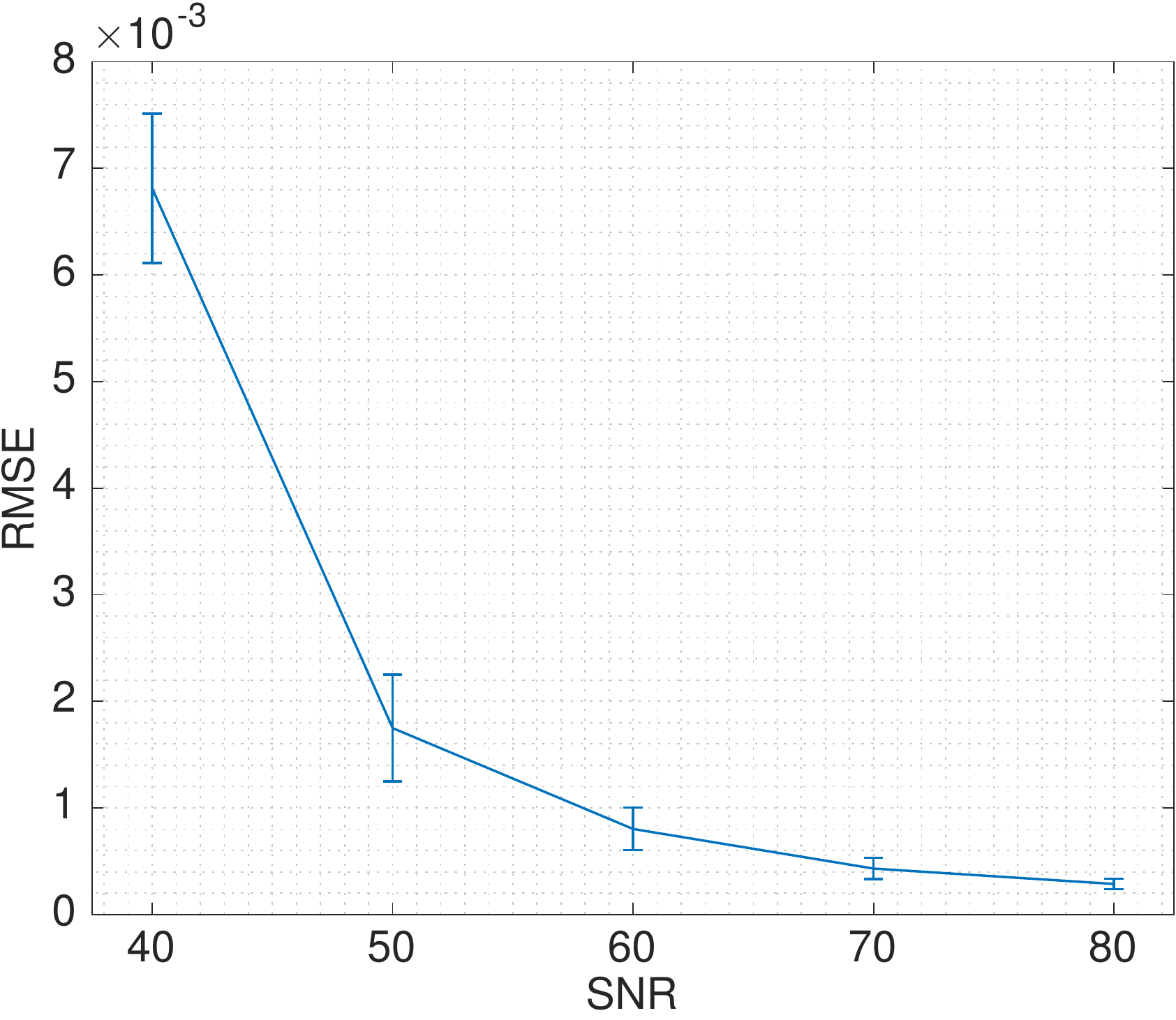} 
	\caption{RMSE for the system state with non-phasor data for $T=2$.}
	\label{fig:rmse}
\end{figure}

We next evaluated how the estimation accuracy of system states depends on the SNR. For this test, the SNRs for non-metered loads and probing data were identical. The state estimation accuracy was evaluated in terms of the root mean square error (RMSE) defined as 
$\sqrt{\sum_{t=1}^T\|\bv_t-\hbv_t\|_2^2/(NT)}$
averaged over 20 Monte Carlo tests. Figure~\ref{fig:rmse} shows how the RMSE decreases for increasing SNR. 

\begin{figure}[t]
	\centering
	\includegraphics[scale=0.4]{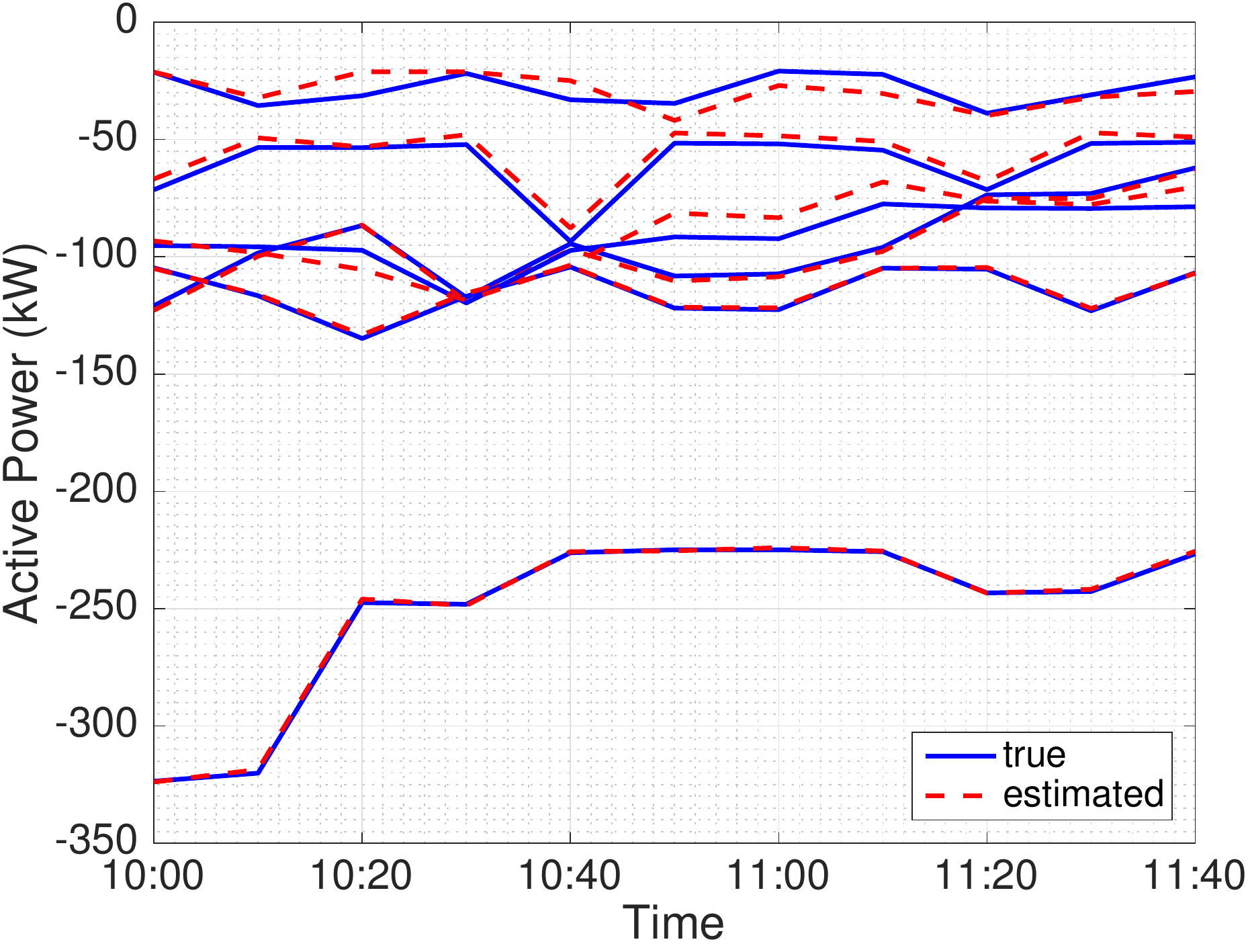}\\ \vspace*{1em}
	\includegraphics[scale=0.4]{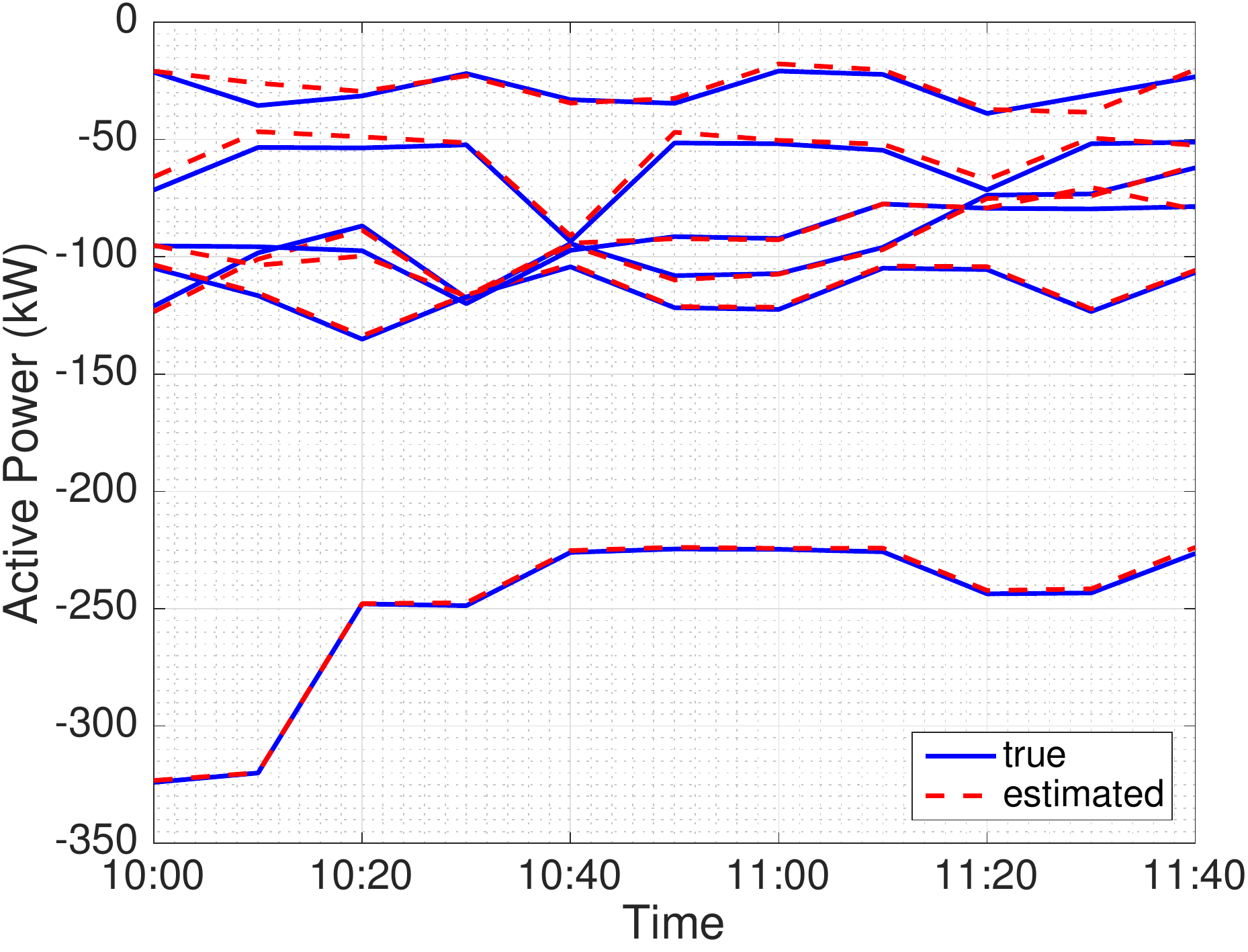}
	\caption{Active power injection estimates using probing with phasor data for $T=4$ (top) and $T=6$ (bottom).}
	\label{fig:PQ4}
\end{figure}

\begin{figure}[t]
	\centering
	\includegraphics[scale=0.4]{PMUn-T6O6P.pdf}\\ \vspace*{1em}
	\includegraphics[scale=0.4]{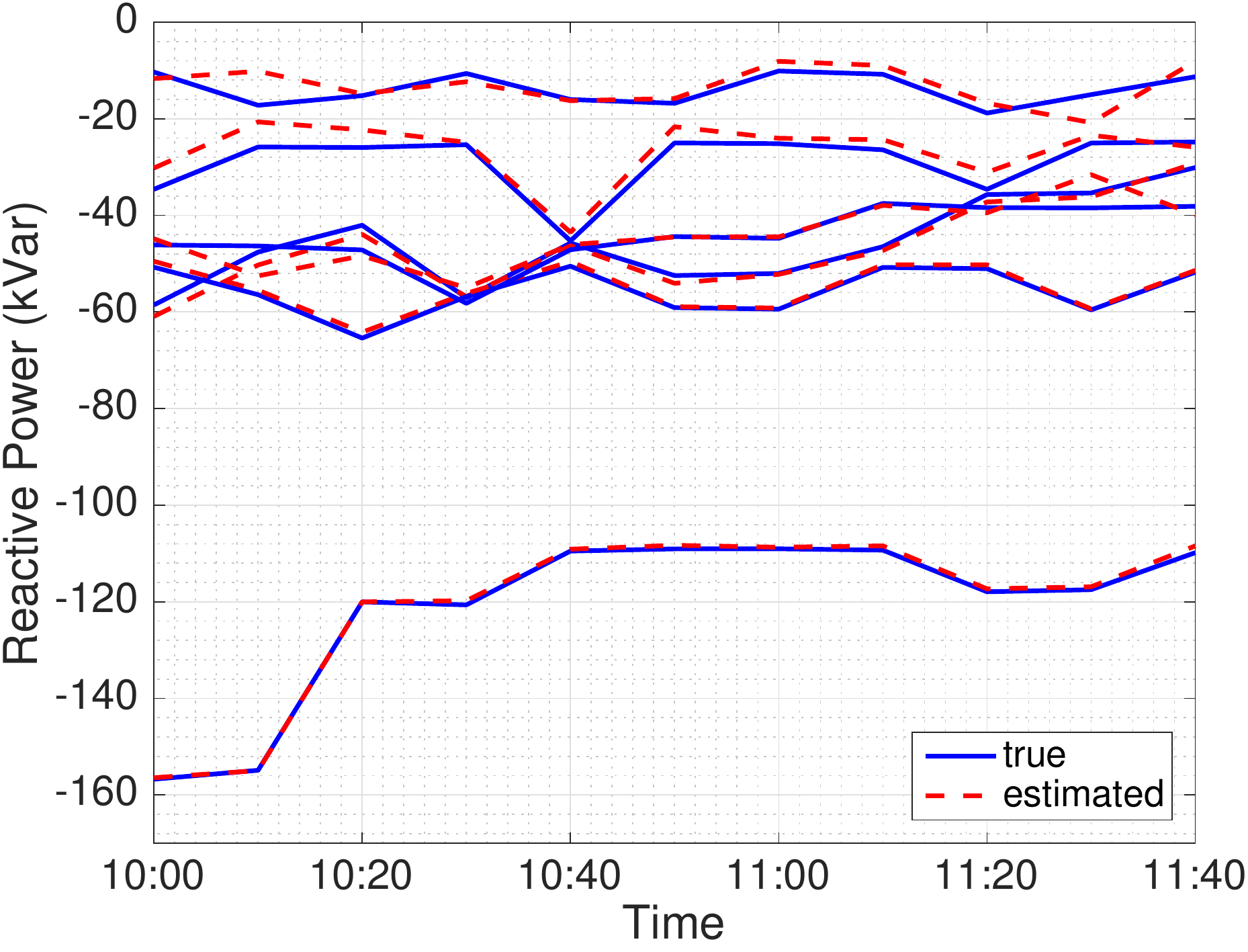}
	\caption{Active (top) and reactive (bottom) power injection estimates using probing with phasor data for $T=6$ and $O=6$.}
	\label{fig:PQ5}
\end{figure}

To validate P2L over different loading conditions, we ran numerical tests for the period of 10:00 a.m. and 01:40 p.m. and every 10 min using phasor data. The SNRs for probing data and non-metered loads were again fixed to $80$ and $60$dB, respectively. Figures~\ref{fig:PQ4} and \ref{fig:PQ5} present the actual and estimated non-metered (re)active loads on buses $\{4,6,15,20,27,31\}$ for $T=4$ and $T=6$, accordingly. The plots show the load estimation improvement by increasing $T$.

Regarding the runtime of our algorithms, each P2L task took between $95-180$~sec, which were allocated as follows:
\begin{itemize}
	\item The linear programs of \eqref{eq:farkas} took $70$~sec overall, to check the feeder compliance of $K=100$ candidate probing setpoints.
	\item The quadratic program of \eqref{eq:MSDr} needed to select the $T$ most diversifying setpoints was solved in less than $10$~sec. 
	\item The SDP formulation of \eqref{eq:CPSSE} together with the heuristic of \cite{MALB16} to obtain a rank-one solution took $25-100$~sec. 
\end{itemize}
The load learning task for the single-slot probing setup $(T=1)$ was solved in less than $15$~sec.

\section{Conclusions}\label{sec:conclusion}
The novel data acquisition scheme of probing an electric grid via smart inverters to infer non-metered loads has been presented. Part I studied the topological observability of grid probing using (non)-phasor data in potentially meshed networks. If a probing setup is deemed topologically observable, Part II has presented a systematic methodology for designing probing injections. The goal is improved estimation accuracy and adherence to inverter and feeder constraints even without knowing non-metered loads. The computational tasks involved in grid probing have been cast as penalized SDP-based solvers and account for noisy measurements and non-stationary loads. 

Numerical tests using synthetic and real-world data on benchmark feeders demonstrate the ensuing take-away simulation findings: \emph{i)} High-accuracy phasor data are better for load recovery than non-phasor data; \emph{ii)} Having the most diverse system states during probing yields better load estimates; \emph{iii)} Probing seemed to yield better estimates under broad voltage regulation range and tight load uncertainty \emph{iv)} Although increasing $T$ improved the system state accuracy, the obtained load estimates were not always better, especially for larger $O$. Nevertheless, we were able to recover a reasonable number of loads; and \emph{v)} Including the extra constraints to strengthen the SDP relaxation provided better numerical accuracy. 

Several questions remain open. Developing scalable solvers perhaps along the lines of \cite{Gang18}; incorporating measurement from distribution lines and transformers~\cite{SevlianR15}; and applying our topological observability framework to detect data attacks in distribution grids; all constitute pertinent research directions. 

\balance
\bibliographystyle{IEEEtran}
\bibliography{myabrv,power}

\end{document}